\documentclass[preprint,12pt,twoside,a4paper]{article}
\usepackage[english]{babel}
\usepackage[utf8x]{inputenc}
\usepackage{array}
\usepackage[a4paper,top=3cm,bottom=2cm,left=3cm,right=3cm,marginparwidth=2cm]{geometry}
\usepackage{fancyhdr} 
\usepackage{amsthm,amssymb,xcolor,amsmath,amsfonts,graphicx,multicol,nccmath}
\usepackage[colorinlistoftodos]{todonotes}
\usepackage[colorlinks=true, allcolors=blue]{hyperref}
\usepackage[export]{adjustbox}
\usepackage{multirow}
\usepackage{caption}
\usepackage{subcaption}
\usepackage{float}
\newtheorem{assumption}{Assumption}
\newtheorem{remark}{Remark}
\newtheorem{proposition}{Proposition}
\newtheorem{theorem}{Theorem}[section]
\newtheorem{definition}{Definition}[section]

\numberwithin{equation}{section}
\title{Jump-induced mixed-mode oscillations through picewise-affine maps	}
\author{Yiorgos Patsios			\\
		Hasselt Universiteit		\\
		\and 
		Renato Huzak				\\
		Hasselt Universiteit		\\	
		\and 
		Nikola Popovi\'c			\\
		University of Edinburgh	\\
		\and
		Peter De Maesschalck		\\
		Hasselt Universiteit		\\	
		}

\date{\today}

\begin{document}

\maketitle
%
%
%
%
%
%

\begin{abstract}
Mixed-mode oscillations (MMOs) are complex oscillatory patterns in which large-amplitude relaxation oscillations (LAOs) alternate with small-amplitude oscillations (SAOs). 
MMOs are found in singularly perturbed systems of ordinary differential equations of slow-fast type, and are typically related to the presence of so-called folded singularities and the corresponding canard trajectories in such systems. 
Here, we introduce a canonical family of three-dimensional slow-fast systems that exhibit MMOs which are induced by relaxation-type dynamics, and which are hence based on a ``jump mechanism", rather than on a more standard canard mechanism. 
In particular, we establish a correspondence between that family and a class of associated one-dimensional piecewise affine maps (PAMs) which exhibit MMOs with the same signature. 
Finally, we give a preliminary classification of admissible mixed-mode signatures, and we illustrate our findings with numerical examples.
\end{abstract}

\section{Introduction}	\label{intro}
In the theory of dynamical systems, one is generally interested in the qualitative behaviour of solutions of differential equations.
Thus, for instance, one investigates bifurcations of equilibria and periodic orbits in dependence of parameters in these systems.
In this paper, we focus on singularly perturbed three-dimensional systems of ``slow-fast" type, with two slow variables and one fast variable.
Such systems are characterised by the variables evolving on different time-scales, which can, in some circumstances, give rise to canard phenomena.
Canards \cite{BENOIT} arise when trajectories of a singularly perturbed system follow an attracting manifold, pass through a folded singularity and then -- somewhat counterintuitively -- stay close to a repelling slow manifold for some time. 
In planar slow-fast systems, the canard phenomenon is often linked to the presence of a (singular) Hopf bifurcation at a turning (fold) point; one typical example is given by the singularly perturbed van der Pol equation \cite{DGKKOW12,DR96}. 
Canards have been studied extensively over the past decades; their study has mainly been based on non-standard analysis \cite{BCDD81,D84, D94}, matched asymptotic expansions, and a geometric approach that combines Fenichel's geometric singular perturbation theory (GSPT) \cite{F72,F79}
and the so-called blow-up technique, which was introduced in the pioneering works of Dumortier and Roussarie \cite {DR96}, as well as of Krupa and Szmolyan \cite {KS01}.
\smallskip

In three-dimensional slow-fast systems with two slow variables, the canard phenomenon can give rise to mixed-mode oscillatory dynamics. Mixed-mode oscillations (MMOs) typically consist of 
large-amplitude oscillations (LAOs) of relaxation type, followed by small-amplitude oscillations (SAOs). While no generally accepted rigorous definition of MMOs seems to exist, a clear, intuitive separation between LAOs and SAOS seems to be evident in most cases;
what draws immediate attention is the pattern that emerges in the alternation between oscillations
of distinct amplitudes. Specifically, a (periodic) MMO is said to have {\em signature} $L^s$ if 
the corresponding orbit undergoes $s$ SAOs, followed by $L$ LAOs, at which point 
that sequence repeats. See \cite{DGKKOW12} for a recent review of this complex oscillatory dynamics, as well as \cite{FG11,JMBNR13,chemists,MN09,MSLG97,VL97} for a small selection of biological, chemical, and physical models in which a variety of MMO patterns have been observed.
\medskip

Among the numerous mechanisms that have been proposed to explain mixed-mode oscillatory dynamics in singularly perturbed systems of slow-fast type, the canard-based mechanism \cite{BKW06,KPK08,W05} has been among the most popular.
Roughly speaking, it combines local passage through the vicinity of a folded singularity -- which explains the SAO component of the corresponding MMO -- with a global return mechanism which results in relaxation (LAO), returning the flow to the basin of attraction of the folded singularity \cite{DGKKOW12}. In the present paper, we introduce an alternative mechanism for the generation of mixed-mode dynamics in three-dimensional slow-fast systems, which we will refer to as the ``jump mechanism". In the process, we will show that the occurrence of MMOs is not necessarily caused by the presence of a folded singularity, as in the canard-based mechanism; in fact, the main characteristic of the MMOs studied in this paper is that 
both LAOs and SAOs are now of relaxation type and that the amplitude of the latter is thus of order $O(1)$ in the singular perturbation parameter.
\smallskip

Our study is inspired by previous work of Szmolyan and Wechselberger \cite{SW04}, Krupa, Popovi\'c, and Kopell \cite{KPK08}, and Rajpathak, Pillai, and Bandyopahdyay \cite{RPB12}. By considering a prototypical family of slow-fast systems which incorporates two jump mechanisms of the type studied in \cite{SW04}, we reproduce MMOs that alternate between LAOs and SAOs of relaxation type; see Figure~\ref{fig:assumptions} for an illustration of the resulting geometry, as well as Section~\ref{sec:assumptions} for a precise definition of our family. As we rely on established results from \cite{SW04}, we do not explicitly need to perform a family blow-up in order to desingularise the flow near fold curves along which normal hyperbolicity is lost. 
In the process, we reduce the study of mixed-mode dynamics in our prototypical family to that of one-dimensional piecewise affine maps (PAMs) \cite{BKYY2000,DRBA08,PRV16,YO96}. In particular, we show that the singular limit of the corresponding first return (or Poincar\'e) map yields a PAM; see Proposition \ref{return-limit-prop}. Piecewise maps  \cite{DIBERNARDO} have been popularised in the study of dynamical systems in recent decades, with a particular focus on models for switching phenomena such as electrical circuits \cite{H00,IK2020,SSI15} and neurons \cite{IMMTDZ11,JMBNR13, GKE17}; 
such maps are naturally related to the corresponding Poincar\'e maps in oscillatory systems.
\smallskip

Thus, we establish a natural two-way correspondence between the family of three-dimensional slow-fast systems studied here and a suitably defined class of one-dimensional PAMs which is associated to the reduced flow on the critical manifold of that family. Specifically, we show that a slow-fast system which satisfies the assumptions in Section~\ref{sec:assumptions} exhibits a periodic MMO of a given signature $L^s$ if the PAM which is associated to that system exhibits a periodic 
MMO with the same signature (Theorem \ref{association_thm}). 
Conversely, we show that any PAM within a very broad class of maps is associated to a slow-fast system within the family defined in Section~\ref{sec:assumptions} (Theorem~\ref{existence_thm}). 
\smallskip

This paper is organised as follows. In Section~\ref{sec:assumptions}, we define the 
three-dimensional family of slow-fast systems which underlies our results; in particular, we 
introduce the jump-type mechanism that gives rise to mixed-mode dynamics in our context.
We state our main results in Section~\ref{statement of the results}; before proving those in 
Section~\ref{section-proofs-all}, we explain how to compute the associated PAM (Section~\ref{sec:derivePAM}). 
In Section~\ref{sec:numerics}, we consider a particular representative from our family of systems to verify our findings numerically, and we present relevant simulations. Moreover, we give numerical evidence of mixed, ``crossover" signatures.
Finally, we end our work with a concluding discussion of our findings, as well as with an outlook to potential future research endeavours.


\section{Slow-fast model and assumptions}\label{sec:assumptions}

We study MMOs in the context of the following canonical three-dimensional family of slow-fast systems in the standard form of geometric singular perturbation theory,
\begin{subequations}\label{VFfast}
\begin{align}
x' &=y-F(x,z,\epsilon,\delta)=:f(x,y,z,\epsilon,\delta), \label{VFfast-a} \\
y' &=\epsilon g_1(x,y,z,\epsilon,\delta), \label{VFfast-b} \\
z' &= \epsilon g_2(x,y,z,\epsilon,\delta); \label{VFfast-c}
\end{align}
\end{subequations}
here, $f$, $g_1$, $g_2$  are $C^{\infty}$-smooth functions in their arguments that will be specified in the following and $\epsilon \geq 0$ is a (small) singular
perturbation parameter. 
We emphasise that \eqref{VFfast} contains an additional parameter $\delta$, the relevance of which will become evident below. 
Correspondingly, $x \in \mathbb{R}$ is a fast variable, while $(y,z) \in \mathbb{R}^{2}$ are slow variables, all of which depend on the fast time $t$. 
To avoid unnecessary abstraction, we have assumed that $f$ is of the specific form $f=y-F$ in \eqref{VFfast-a}; that assumption is made without loss of generality and does not represent a major restriction. Our study
will chiefly be based on Fenichel's geometric singular perturbation theory (GSPT) \cite{F72,F79}; an excellent introduction can be found in \cite{JONES}. 

Rewriting the above
{\it fast system} in terms of the slow time variable $\tau=\epsilon t$, we obtain the equivalent 
{\it slow system}
\begin{subequations}\label{VFslow}
\begin{align}
\epsilon\dot{x} &=y-F(x,z,\epsilon,\delta), \\
\dot{y} &=g_1(x,y,z,\epsilon,\delta), \\
\dot{z} &=g_2(x,y,z,\epsilon,\delta),
\end{align}
\end{subequations}
where the overdot denotes differentiation with respect to $\tau$.

In the singular limit of $\epsilon=0$, the above systems yield the {\it layer problem}
\begin{subequations}\label{VFlayer}
\begin{align}
x' &=y-F(x,z,0,\delta), \\
y' &=0, \\
z' &=0
\end{align}
\end{subequations}
and the {\it reduced problem}
\begin{subequations}\label{VFreduced}
\begin{align}
0 &=y-F(x,z,0,\delta), \\
\dot{y} &=g_1(x,y,z,0,\delta), \\
\dot{z} &=g_2(x,y,z,0,\delta),
\end{align}
\end{subequations}
respectively. In particular, \eqref{VFreduced} allows us to define the ($\delta$-family of) critical manifolds $\mathcal{S}:= \lbrace (x,y,z) \in \mathbb{R}^{3}: f(x,y,z,0,\delta)=0 \rbrace $, which is of central importance in GSPT: the sign of $\frac{\partial f}{\partial x}$ determines the stability of the steady states of the layer problem in \eqref{VFlayer}, which are located on $\mathcal{S}$.
Specifically, orbits that are initiated away from $\mathcal{S}$ will converge to attracting branches of the critical manifold under the layer flow of \eqref{VFlayer}; on $\mathcal{S}$,
they will then be subject to the reduced flow of \eqref{VFreduced}.
Away from zeros of $\frac{\partial f}{\partial x}$, the critical manifold $\mathcal{S}$ is {\it normally hyperbolic}; 
by Fenichel's First Theorem, normally hyperbolic segments of $\mathcal{S}$ will perturb, for $\epsilon$ positive and sufficiently small, to a slow manifold $\mathcal{S}_\epsilon$ \cite{F79}. 
Correspondingly, the reduced flow on $\mathcal{S}$ will perturb in a regular fashion to the slow flow on $\mathcal{S}_\epsilon$. 
Likewise, the fast flow of \eqref{VFfast} will be a regular perturbation of the layer flow off $\mathcal{S}$.

A canonical scenario in which normal hyperbolicity is lost is found at so-called fold curves in \eqref{VFfast}, where orbits exhibit jumping behaviour. 
The reduced flow on $\mathcal{S}$ is directed towards (attraction) or away from (repulsion) these fold curves, which results in orbits having to jump to a different segment of the critical manifold there. 
When such behaviour occurs in a periodic fashion, {\it relaxation oscillation} is observed. 
The emergence of fold-induced relaxation oscillation in three-dimensional slow-fast systems was studied in detail in \cite{SW04}, where the desingularisation technique known as ``blow-up" \cite{DR96} was applied to remedy the loss of normal hyperbolicity.

We follow the same approach here and proceed to make the following (analogous) assumptions on the singular geometry of Equation~\eqref{VFfast}.

\begin{figure}[H]
\includegraphics[width=14cm,center]{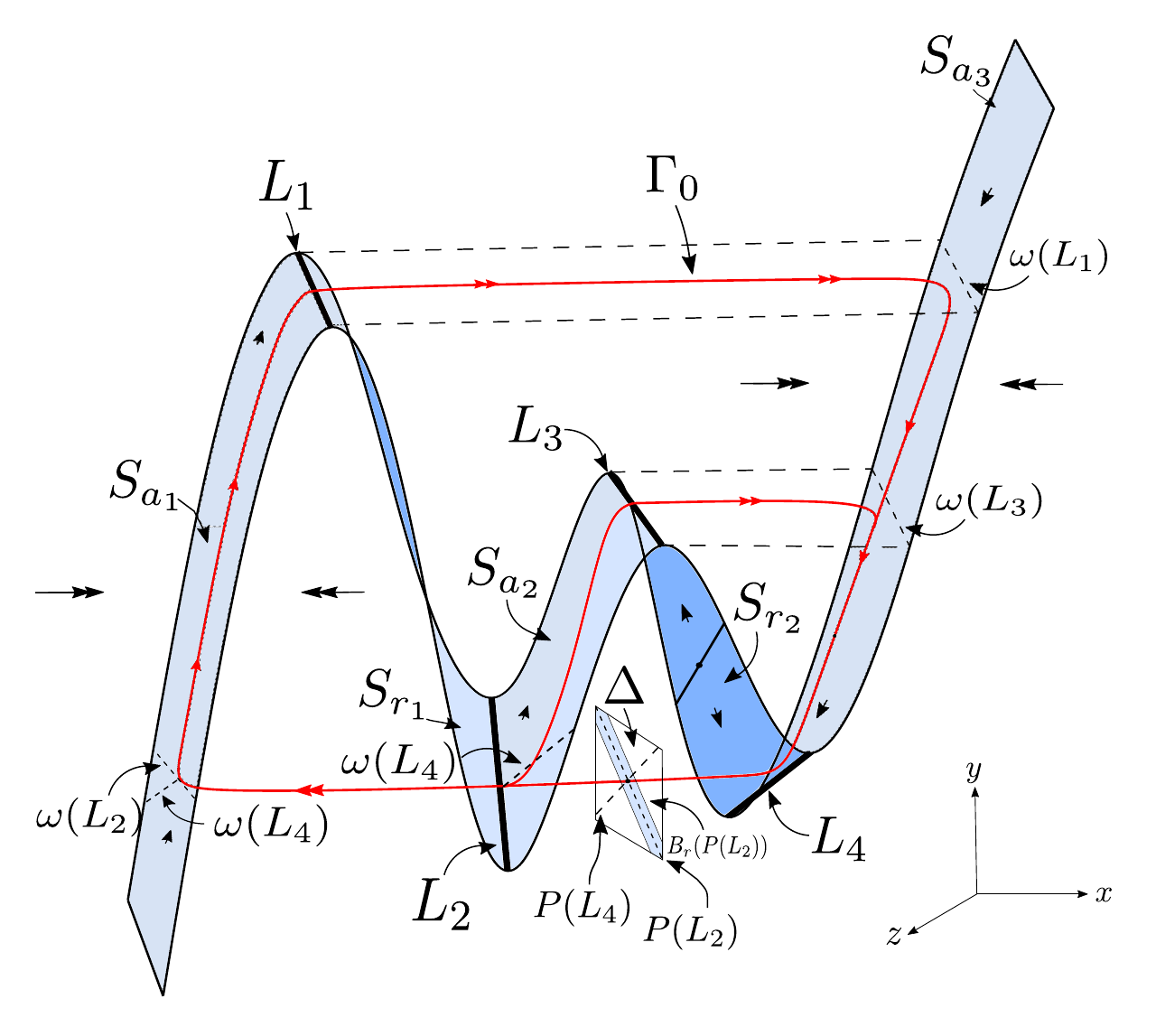} 
\caption{Critical manifold $\mathcal{S}$ with ``ambiguous'' singular orbit $\Gamma_0 = \Gamma_0^L\cup\Gamma_0^S$ and section $\Delta$.} 
\label{fig:assumptions}
\end{figure}

\begin{assumption}\label{A1} 	
Possibly restricting to a part of phase space for $z\in (z_{\rm min},z_{\rm max})$, we assume that the critical manifold $\mathcal{S}$ is $\xi$-shaped, or ``Bactrian"-shaped; see Figure~\ref{fig:assumptions}. In other words, $\mathcal{S}$ can be written as
\begin{align*}
\mathcal{S} = \mathcal S_{a_{1}} \cup L_{1} \cup \mathcal S_{r_{1}} \cup L_{2} \cup \mathcal S_{a_{2}} \cup L_{3} \cup \mathcal S_{r_{2}} \cup L_{4} \cup \mathcal S_{a_{3}},
\end{align*}
where $\mathcal S_{a_{1}} \cup \mathcal S_{a_{2}} \cup \mathcal S_{a_{3}} = \mathcal{S}\cap \lbrace \frac{\partial f}{\partial x}(x,y,z,0, \delta)<0 \rbrace $
and $\mathcal S_{r_{1}} \cup \mathcal S_{r_{2}} = \mathcal{S} \cap \lbrace \frac{\partial f}{\partial x}(x,y,z,0,\delta)>0 \rbrace $ denote the normally attracting and normally repelling segments of
$\mathcal S$, respectively, which are divided by four ($\delta$-families of) fold curves along which normal hyperbolicity is lost, denoted $L_{1}$, $L_{2}$, $L_{3}$ and $L_{4}$ from left to right. These fold curves can be written as graphs
\begin{align*}
L_{i} := \lbrace (x,y,z)\in \mathcal{S}: (x,y,z)=(\nu _{i}(z,\delta),\phi(\nu _{i}(z,\delta),z),z)\rbrace\quad  \text{for }i=1,2,3,4;
\end{align*}
here, $\phi$ and $\nu_{i}$ ($i=1,2,3,4$) are appropriately defined functions along which the non-degeneracy conditions
\begin{align*}
\frac{\partial f}{\partial x}\big(\nu _{i}(z,\delta),\phi(\nu _{i}(z,\delta),z),z, 0, \delta\big)=0\quad\text{and}\quad
\frac{\partial^2 f}{\partial x^2}\big(\nu _{i}(z,\delta),\phi(\nu _{i}(z,\delta),z),z, 0, \delta\big)\neq 0
\end{align*}
are satisfied.

As $\nu _{1}<\nu _{2}<\nu _{3}<\nu _{4}$, by assumption, the normally attracting segments of $\mathcal S$ can equally be represented as
\begin{align*}
\mathcal S_{a_{1}} &:= \lbrace (x,y,z)\in \mathcal{S} \cap \lbrace x<\nu _{1} \rbrace \rbrace,\quad
\mathcal S_{a_{2}} := \lbrace (x,y,z)\in \mathcal{S} \cap \lbrace \nu _{2}<x<\nu _{3} \rbrace \rbrace,\quad\text{and} \\
\mathcal S_{a_{3}} &:= \lbrace (x,y,z)\in \mathcal{S} \cap \lbrace x>\nu _{4} \rbrace \rbrace,
\end{align*}
while the normally repelling ones are given by
\begin{align*}
\mathcal S_{r_{1}} := \mathcal{S}\cap \lbrace (x,y,z): \nu _{1}<x<\nu _{2} \rbrace\quad\text{and}\quad
\mathcal S_{r_{2}} := \mathcal{S}\cap \lbrace (x,y,z): \nu _{3}<x<\nu _{4}  \rbrace.
\end{align*}
\end{assumption}

\begin{assumption}[Normal switching condition]\label{A2}	
We assume that
\begin{equation}\label{eq:nsc}
f_{y} g_{1} + f_{z} g_{2} \Big\rvert_{p \in L_{i}} \neq 0\qquad\text{for }i=1,2,3,4,
\end{equation}
i.e., that any fold point $p$ on $L_i$ is a jump point. 
In other words, \eqref{eq:nsc} asserts that the reduced flow of \eqref{VFreduced} is unbounded on the fold lines $L_i$ and that orbits must hence jump there.
It is therefore required that the reduced flow on both sides of the fold lines $L_{i}$ ($i=1,3,4$) is transverse to, and directed towards, the fold lines $L_{i}$ at all times. 
\end{assumption}

We define by $\omega(L_{2})$ the projection of the fold line $L_{2}$ onto the attracting sheet $\mathcal S_{a_{1}}$ of $\mathcal{S}$; moreover, we define $\omega(L_{1})$ and $\omega(L_{3})$ as the projections of the fold lines $L_{1}$ and $L_{3}$ onto $\mathcal S_{a_{3}}$. Likewise, we define the projection $\omega(L_{4})$ of the fold line $L_{4}$ onto both $\mathcal{S}_{a_{1}}$ and $\mathcal{S}_{a_{2}}$. (In spite of the suggestive notation, these projections should not be confused with $\omega$-limit sets of these fold lines.) Then, we assume the following.
\begin{assumption}[Transversality of reduced flow]\label{A4}
For $\delta=0$, the reduced flow of \eqref{VFreduced} is transverse to $\omega(L_{1})$ and $\omega(L_{3})$ on $\mathcal S_{a_{3}}$, transverse to $\omega(L_{2})$ on $\mathcal S_{a_{1}}$, and transverse to $\omega(L_{4})$ on both $\mathcal{S}_{a_{1}}$ and $\mathcal{S}_{a_{2}}$.
\end{assumption}

Let $\Delta$ denote a section between $L_{2}$ and $L_{4}$ that is transverse to the layer flow of \eqref{VFlayer}, and let $P(L_{2})$ and $P(L_{4})$ denote the projections of the fold lines $L_{2}$ and $L_{4}$, respectively, onto $\Delta$.
The two projection lines intersect transversally at a point $P_{c} := P(L_{2}) \cap P({L_{4}})$, as indicated in Figure~\ref{fig:assumptions}.
Then, we make the following assumption:
\begin{assumption}[Breaking mechanism]\label{A3}	
For {$\delta>0$} sufficiently small and $\epsilon=0$, the $z$-parametrized curves $P(L_{2})$ and $P(L_{4})$ intersect transversely at some $z$-value $z_0(\delta)$, with $z_0(0)=z_0$ and $z'_0(0)=0$; in particular, the intersection point $P_c$ between the curves thus depends on $\delta$.
Furthermore, we assume that for $z < z_{0}(\delta)$, the evolution of $P(L_{4})$ in forward time remains below the normally hyperbolic sheets $\mathcal S_{r_{1}}$ and $\mathcal S_{a_{2}}$, and that it lands directly on the opposite attracting sheet $\mathcal S_{a_{1}}$.
On the other hand, for $z>z_{0}(\delta)$, the evolution of $P(L_{4})$ lands on $\mathcal S_{a_{2}}$, safely away from the fold line $L_{2}$; see again Figure~\ref{fig:assumptions}.
\end{assumption}
We assume that $z_0'(0)=0$ to ensure that the piecewise affine map (PAM) associated to \eqref{VFfast}, as introduced in Section~\ref{statement of the results}, has a jump at the origin. As will become clear there, that assumption is made without loss of generality: the general case yields a piecewise affine map with a jump at non-zero $z_0'(0)$, which can be studied in an analogous fashion after a translation. 


Under the above assumptions, we can already give a partial slow-fast analysis of the system in \eqref{VFfast}; it will be convenient to do so now in order to explain our final assumption. Consider first the fold curve $L_4$ where we assume all orbits to jump. Given Assumption~\ref{A3}, the $z$-coordinate of a given orbit will determine whether it is attracted to $\mathcal{S}_{a_1}$ or to  $\mathcal{S}_{a_2}$ under the layer flow. For $z=z_0$, the fate of orbits cannot be decided; consideration of the perturbation terms in \eqref{VFfast} and a blow-up of $L_2$ would be necessary to describe the flow in that case. We expect that in general, canard phenomena are possible near $z=z_0$ where orbits will follow part of $\mathcal{S}_{r_1}$; here, we do not consider that scenario. However, as is clear from the above discussion, the $z$-value $z_0$ will nevertheless play a central role in our analysis. We will highlight two possible singular orbits passing through the transcritical intersection point $P_c=P(L_2)\cap P(L_4)$ at $z=z_0$: one orbit will continue along $\mathcal{S}_{a_1}$, while the other will continue along $\mathcal{S}_{a_2}$. In both instances, we will assume that the sought-after singular orbit follows the reduced flow until a fold line is reached -- $L_1$ in the former case and $L_3$ in the latter -- at which point the orbit jumps to $\mathcal{S}_{a_3}$ and ultimately reaches $L_4$, following again the reduced flow. 
This ``ambiguous" behaviour of the singular flow is a key point in our study and motivates the following Assumption~\ref{A5}.

\begin{assumption}[Ambiguous singular orbit]\label{A5}
There exists a singular closed orbit $\Gamma_0^L$ that is defined by concatenating the reduced flow on $\mathcal{S}_{a_1}$ and on $\mathcal{S}_{a_3}$ with the layer flow between $L_1$ and $\mathcal{S}_{a_3}$ and between $L_4$ and $\mathcal{S}_{a_1}$, respectively.
Further, there exists a singular closed orbit $\Gamma_0^S$ that is defined by concatenating the reduced flow on $\mathcal{S}_{a_2}$ and on $\mathcal{S}_{a_3}$ with the layer flow between $L_3$ and $\mathcal{S}_{a_3}$ and between $L_4$ and $\mathcal{S}_{a_2}$, respectively. 
Both $\Gamma_0^L$  and $\Gamma_0^S$ contain the point of intersection $P_c$ defined in Assumption~\ref{A3}; see Figure~\ref{fig:assumptions}. 
Finally, we define $\Gamma_0 = \Gamma_0^L\cup\Gamma_0^S$ as the ``ambiguous'' singular orbit.
\end{assumption}
\begin{remark}
Any orbit $\Gamma_{0}$ in the above family forms a natural boundary between oscillations of different amplitude.
\end{remark}
For the sake of convenience, we will assume that both singular orbits $\Gamma_0^L$ and $\Gamma_0^S$ lie in a plane $z=z_0$. It follows that $\dot{z}|_{(z,\delta)=(z_0,0)}=0$ and, hence, that we can write
\begin{equation}\label{eq:g2}
g_2(x,y,z,\epsilon,\delta) = \delta G(x,y,z,\epsilon,\delta) + (z-z_0)H(x,y,z,\epsilon,\delta)
\end{equation}
in \eqref{VFfast-c}.
Our aim in this paper is to study the behaviour of orbits near $\Gamma_0$. To that end, we will formulate a first return map on the section $\Delta$ defined above which is transverse to the layer flow.  We will present a partial study here: to be precise, we will restrict to characterising orbits that are sufficiently close, but not too close, to the point $P_c = P(L_2)\cap P(L_4)$; in other words, we will consider orbits in a sufficiently small neighbourhood of $P_c$ inside $\Delta$, uniformly away from $P(L_2)$. Our analysis will rely in part on \cite{SW04}, which will allow us to describe the persistence of both $\Gamma_0^L$ and $\Gamma_0^S$.

\section{Statement of results}\label{statement of the results}


Recall the definition of the section $\Delta$ which is located between the fold lines $L_2$ and $L_4$ and which is transverse to the layer flow of Equation~\eqref{VFlayer}.  
Also, recall that the projection of the fold lines $L_2$ and $L_4$ onto $\Delta$ is denoted by $P(L_2)$ and $P(L_4)$, respectively.   
Our first result concerns the well-definedness of the first return map from $\Delta$ to itself under the flow of Equation~\eqref{VFfast} and is relatively straightforward, 
since most of the relevant dynamics occurs along hyperbolically attracting parts of the critical manifold $\mathcal{S}$. 
We do, however, need to take additional care in a neighbourhood of $P(L_2)$ on $\Delta$, as the fate of orbits sufficiently close to the fold curve $L_2$ is difficult to analyse: such orbits could either jump onto an attracting sheet or follow a repelling sheet of the critical manifold $\mathcal{S}$ after passing near $L_2$, which would give rise to canard behaviour.  
As stated above, we will avoid this unpredictability here; we will therefore formulate a result on the first return map that avoids a neighbourhood of $P(L_2)$.  

First, we note that the form of the vector field in \eqref{VFfast} allows us to conclude that $P(L_2)$ and $P(L_4)$ are smooth $\delta$-families of graphs
\begin{align*}
y= \phi_{L_2}(z,\delta)\qquad\text{and}\qquad y = \phi_{L_4}(z,\delta)
\end{align*}
which intersect in the point $P_c$ given by $(y_0,z_0,0) + O(\delta)$.
Next, for $r>0$, we define the open neighbourhood
\begin{align*}
B_r(P(L_2)) = \{ (y,z)\in\Delta : |y-\phi_{L_2}(z,\delta)|  < r\}
\end{align*}
of $P(L_2)$. Then, we have the following result.
    
\begin{theorem}\label{firstreturnmap_thm} 
There exists an open neighborhood $\mathcal{U}$ of the point $P_c$ in $\Delta$ such that,
for each $\delta>0$ and $r>0$ sufficiently small, there exists $\epsilon>0$ small enough such that the first return map
\begin{align*}
\Pi: \mathcal{U}\setminus B_{r}(P(L_2))\subset \Delta\to\Delta\colon (y,z) \mapsto (\mathcal{Y}(y,z,\epsilon,\delta),\mathcal{Z}(y,z,\epsilon,\delta))
\end{align*}
is well-defined. Here, we write
\begin{align*}
\mathcal{Y}(y,z,\epsilon,\delta) = \phi_{L_4}(\mathcal{Z}(y,z,0,\delta),\delta) + \mathcal{E}(y,z,\epsilon,\delta)
\end{align*}
for some function $\mathcal{E}$ that is uniformly $o(1)$ as $\epsilon\to 0$. 
\end{theorem}
Theorem~\ref{firstreturnmap_thm} will be proved in Section~\ref{proof of first return map theorem}. Let us now give some heuristics on how the return map $\Pi$ can be related to a suitably defined PAM, which we require in order to formulate our next result. We will give full proofs in Sections~\ref{sec:derivePAM} and \ref{section-proofs-all} below.

In the singular limit of $\epsilon=0$, the map $\Pi$ is given by
\begin{align*}
\Pi|_{\epsilon=0}: (y,z)\mapsto ( \phi_{L_4}(\mathcal{Z}(y,z,0,\delta),\delta),\mathcal{Z}(y,z,0,\delta) );
\end{align*}
since the image of that map lies on $P(L_4)$, it makes sense to also restrict its domain to the graph $y=\phi_{L_4}(z,\delta)$, reducing it in essence to a one-dimensional map
\begin{align*}
\Pi_0: z\mapsto \mathcal{Z}_{1}(z,\delta) := \mathcal{Z}(\phi_{L_4}(z,\delta),z,0,\delta).
\end{align*}
The map $\Pi_0$ is only defined for values of $z$ that are at least an $O(r)$-distance away from $z_0(\delta)$; however, since we can apply Theorem~\ref{firstreturnmap_thm} for any choice of $r$, the $\epsilon=0$-limit of $\Pi$ is actually defined for all $ z \not= z_0(\delta) $.
Let us now consider the neighbourhood of $z=z_0$ by writing $z=z_0 + \delta Z$. Recalling that $z_0\mapsto z_0$ for $\delta=0$, which marks the transverse intersection point of $P_c$ of $P(L_2)$ and $P(L_4)$ in Assumption~\ref{A3}, we arrive at
\begin{align*}
\tilde{\Pi}_0 : Z\mapsto \tilde{\mathcal{Z}}_1(Z,\delta) := \frac{\mathcal{Z}_{1}(z_0+\delta Z,\delta) - \mathcal{Z}_{1}(z_0,0)}{\delta}\qquad\text{for }
z_0+\delta Z\ne z_0(\delta),
\end{align*}
as $z_0=\mathcal{Z}_{1}(z_0,0)$.
Using Assumption~\ref{A3} once more, we can assume the map $\tilde{\Pi}_0$ to have a well-defined limit
\[
\tilde{\Pi}_{00} : Z\mapsto \lim_{\delta\to0}\tilde{\mathcal{Z}}_1(Z,\delta)\qquad\text{for }Z\not=0.
\]
Since $O(Z^2)$-terms are scaled away in the above limit, one expects the map $\tilde{\Pi}_{00}$ to be piecewise affine. While the above argument is heuristic, it can be made rigorous by relation to the vector field in \eqref{VFfast} and on the basis of the two limiting systems that are obtained therefrom for $\epsilon=0$. Below, we express this correspondence in terms of a definition; a rigorous proof of our heuristics can be found in Proposition \ref{return-limit-prop}.

\begin{definition} \label{assosc_def0}
Let $\gamma$ be a curve on a normally hyperbolic segment of the critical manifold $\mathcal{S}$ for  Equation~\eqref{VFfast} that is parametrized by $x$. Then, we define the affine map $M_\gamma$ associated to $\gamma$ as
\begin{align*}
M_\gamma(Z) = \big(e^{\int_{\gamma} p(x)dx}\big)Z + \left(\int_{\gamma} q(x) e^{\int_{\gamma_{x}} p(x')dx'}dx \right),
\end{align*}
with
\begin{align*}
p(x) := \frac{H(x,y,z,0,0)}{g_1(x,y,z,0,0)} F_{x}(x,z)
\quad\text{and}\quad
q(x) := \frac{G(x,y,z,0,0)}{g_1(x,y,z,0,0)} F_{x}(x,z),
\end{align*}
where we substitute $y$ and $z$ with the $y$-coordinate and the $z$-coordinate of $\gamma$, respectively. (In particular, we note that $y=F(x,z):=F(x,z,0,0)$.)
\end{definition}
We remark that the coefficient of $Z$ in the definition of $M_\gamma$ is strictly positive, as well as that the last integral is defined along a curve $\gamma$ that is parametrized from $x$ until the end of the curve; see \eqref{ODE_soln} below.

Next, we apply Definition~\ref{assosc_def0} in the context of Equation~\eqref{VFfast}, i.e., to the slow portions of the singular orbits $\Gamma_0^S$ and $\Gamma_0^L$, which are given by $\Gamma_0^S\cap\mathcal{S}_{a_2}$, $\Gamma_0^S\cap\mathcal{S}_{a_3}$, $\Gamma_0^L\cap\mathcal{S}_{a_1}$, and $\Gamma_0^L\cap\mathcal{S}_{a_3}$.  
\begin{definition}\label{assosc_def}
We define the piecewise affine map (PAM) 
\begin{equation}
M(Z)=\begin{cases}
M_{\Gamma_0^L\cap \mathcal{S}_{a_3}} \circ M_{\Gamma_0^L\cap \mathcal{S}_{a_1}}(Z)=:M_{1}(Z) & \quad\text{for }Z<0, \\
M_{\Gamma_0^S\cap \mathcal{S}_{a_3}} \circ M_{\Gamma_0^S\cap \mathcal{S}_{a_2}}(Z)=:M_{2}(Z) & \quad\text{for }Z>0,
\end{cases}
\label{PAM}
\end{equation}
and we say that $M$ is \emph{associated} with the vector field in \eqref{VFfast}. 
\end{definition}
We refer to Section~\ref{sec:derivePAM} for specific expressions for $M$ in the context of  \eqref{VFfast}.

\begin{remark}
We note that the $Z$-coefficient of the PAM $M$ is strictly positive, as was the case in Definition~\ref{assosc_def0}.  
However, the image of the map could contain $Z=0$; in fact, the situation where the images $M|_{Z<0}$ and $M|_{Z>0}$ partly overlap represents the most interesting scenario here. 
We will see that the heuristics following the statement of Theorem~\ref{firstreturnmap_thm} can be proved rigorously, and we will show that the limiting map $\tilde{\Pi}_{00}$ is precisely the associated PAM $M$. The fact that the singular limit of the first return map $\Pi$ is not one-to-one causes parts of the remaining analysis to differ from \cite{SW04}.
\end{remark}

The following is our second main result:

\begin{theorem}\label{association_thm}
Given a slow-fast system of the form in \eqref{VFfast} that satisfies Assumptions~\ref{A1} through \ref{A5}, assume that its associated PAM $M$, as defined in \eqref{PAM}, exhibits a stable periodic MMO with signature ${L_{1}}^{s_{1}}{L_{2}}^{s_{2}}...{L_{k}}^{s_{k}}$, for some $k \in \mathbb{N}$. Further, assume that this MMO avoids the discontinuity point at $Z=0$. Then,
 \eqref{VFfast} exhibits a stable MMO with the same signature ${L_{1}}^{s_{1}}{L_{2}}^{s_{2}}...{L_{k}}^{s_{k}}$, for $\epsilon,\delta>0$ sufficiently small.
\end{theorem}

By requiring that $\epsilon$ and $\delta$ be (sufficiently) small in Theorem~\ref{association_thm}, we mean that for every $\delta>0$ small, there exists $\epsilon_0(\delta)>0$ small such that the result is true for every $\epsilon\in(0,\epsilon_0(\delta)]$.
We prove Theorem~\ref{association_thm} in Section~\ref{proof of assosciation theorem}. 

Our theoretical results are complemented by an ``inverse" theorem: any PAM is associated to a suitably chosen slow-fast system of the form in \eqref{VFfast}. In the proof of the following theorem, we will make specific choices for the functions $F$, $g_1$, and $g_2$ therein, which will allow us to obtain convenient expressions for the corresponding vector field.

\begin{theorem}\label{existence_thm}
Assume we are given any piecewise affine map (PAM) of the form 
\begin{align*}
M(Z)= \begin{cases}
a_{11}Z + a_{12}  & \quad\text{for }Z < 0, \\
a_{21}Z + a_{22} & \quad\text{for }Z > 0,
\end{cases}
\end{align*}
with $a_{11}, a_{21}>0$. Then, there exists a slow-fast system of the form in \eqref{VFfast} which satisfies Assumptions~\ref{A1} through \ref{A5} such that the PAM associated to the vector field is given by $M$.
\end{theorem}
Again, Theorem~\ref{existence_thm} will be proved in Section~\ref{Proof of existence theorem}.

\begin{remark}
Theorem~\ref{association_thm} concerns MMOs that are stable both for the PAM $M$ and the associated slow-fast system in \eqref{VFfast}. Given the form of the Poincar\'e map $\Pi$, it then necessarily follows that $a_{11}^L a_{21}^s <1$ in the definition of $M$; see Section~\ref{proof of assosciation theorem} for details. (Here, $L = \sum^{k}_{i=1} L_{i}$ and $s = \sum^{k}_{i=1} s_{i}$, as per the notation of Theorem~\ref{association_thm}.) Clearly, this condition is implied by the somewhat more generic requirement that $a_{11},a_{21}<1$, which is for instance imposed in \cite{RPB12}; cf.~also Section~\ref{sec:atmost+atleast} below.
\end{remark}

\section{Computation of associated PAM}\label{sec:derivePAM}

\begin{figure}[H] \label{fig:critical_manifold_slice}
\includegraphics[width=11cm,scale=1.0,center]{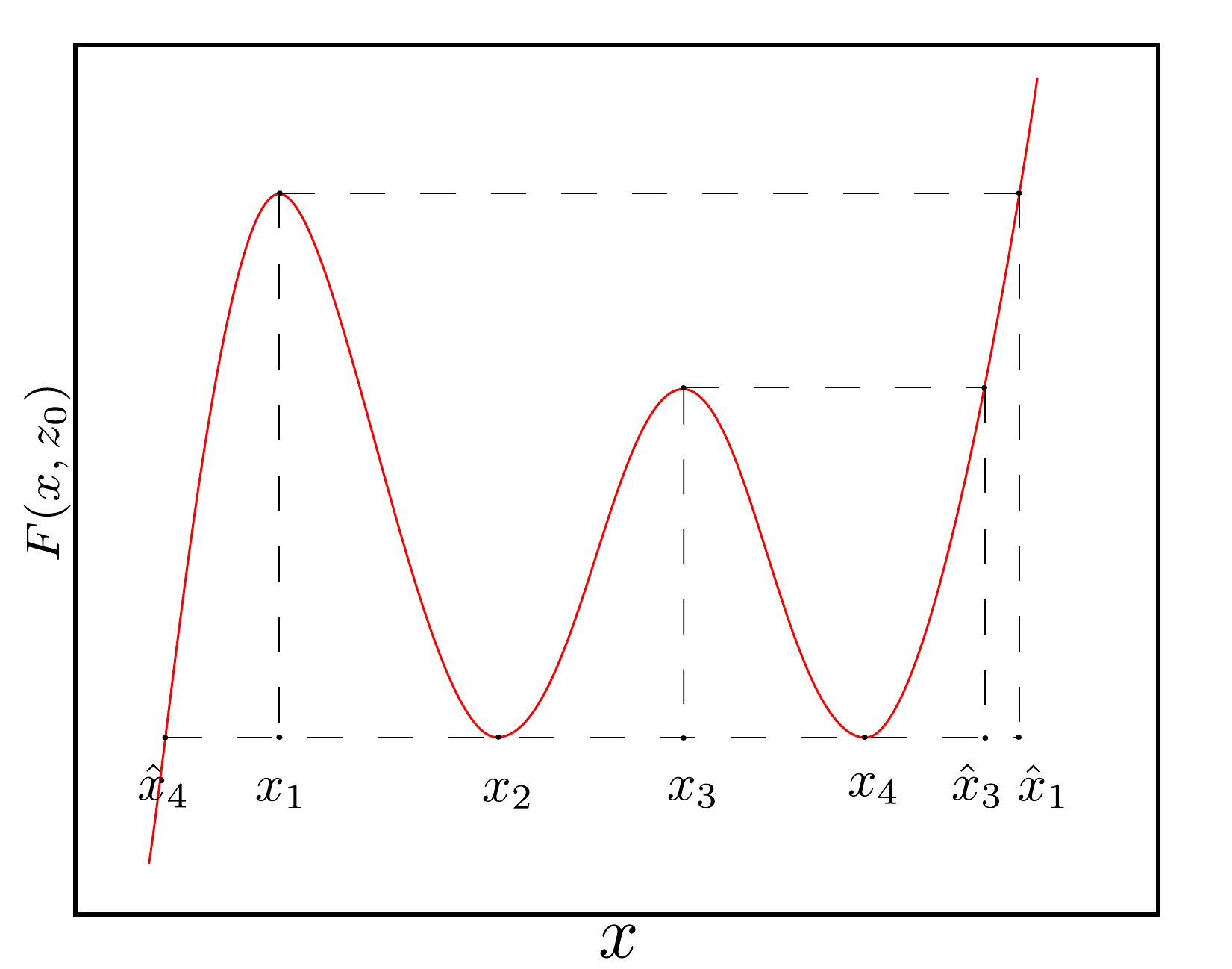}
\caption {The critical manifold $\mathcal{S}$=$\lbrace (x,y,z)\in \mathbb{R}^{3}: y=F(x,z,0,\delta) \rbrace$ for $z=z_{0}$ and $\delta=0$. The function $F(x,z_0):=F(x,z_0,0,0)$ assumes the same value at $\hat{x}_{4}$, $x_{2}$, and $x_{4}$; see Assumption~\ref{A3}. Here, $\hat{x}_{3}$ and $\hat{x}_{1}$ are the $x$-coordinates of the points whose $y$-coordinates are $F(x_{3},z_{0})$ and $F(x_{1},z_{0})$, respectively, which both lie to the right of $x_{4}$.}
\label{fig:critical_manifold}
\end{figure}

In this section, we obtain expressions for the PAM $M$ defined in \eqref{PAM} that is associated to the vector field in \eqref{VFfast}. Moreover, we establish formally that $M$ equals the limit $\tilde{\Pi}_{00}$ of the first return map $\Pi$; see Proposition~\ref{return-limit-prop}. Recall the definition of the system in \eqref{VFfast}, which satisfies Assumptions~\ref{A1} through \ref{A5}. In particular, Assumption~\ref{A1} implies that the function $F(x,z,0,\delta)$ has four distinct local extrema with respect to the variable $x$, which we denote by
\begin{align*}
x_{1}(z,\delta)<x_{2}(z,\delta)<x_{3}(z,\delta)<x_{4}(z,\delta);
\end{align*}
see Figure~\ref{fig:critical_manifold} for an illustration. 


Next, we recall that we assumed for the sake of convenience that $\Gamma_0^S$ and $\Gamma_0^L$ lie in the plane $\{z=z_0\}$; see \eqref{eq:g2} and  the text above that equation. Hence, we can make the general definition of the PAM $M$ somewhat more explicit in the present context. Specifically, the assumption in \eqref{eq:g2} allows us to introduce the following rescaling: since we are interested in $z\approx z_0$ and $\delta$ small in \eqref{VFfast}, we define the transformation
\begin{align}\label{eq:rescale}
z=z_{0} + \delta Z.
\end{align}
Substituting \eqref{eq:rescale} and \eqref{eq:g2} into \eqref{VFfast}, we have
\begin{align*}
x' &=y-F(x,z_{0}+\delta Z,\epsilon,\delta), \\
y' &=\epsilon g_1(x,y,z_{0}+\delta Z,\epsilon,\delta), \\
Z' &=\epsilon\big[ G(x,y,z_0+\delta Z,\epsilon,\delta)+Z H(x,y,z_0+\delta Z,\epsilon,\delta)\big],
\end{align*}
which implies
\begin{align*}
x' &=y-F(x,z_{0})+ \mathcal{O}(\epsilon,\delta), \\
y' &=\epsilon \big[g_1(x,y,z_0,0,0) + \mathcal{O}(\epsilon,\delta)\big], \\
Z' &=\epsilon\big[G(x,y,z_{0},0,0)+ Z H(x,y,z_{0},0,0) + \mathcal{O}(\epsilon,\delta)\big]
\end{align*}
after Taylor expansion of $F$, $g_1$, $G$, and $H$. (Here, we again write $F(x,z_0)=F(x,z_0,0,0)$.)
Reverting to the ``slow time" $\tau$ in the above, we find
\begin{subequations}\label{slow_proof}
\begin{align}				
\epsilon\dot{x} &=y-F(x,z_{0})+ \mathcal{O}(\epsilon,\delta), \\
\dot{y} &= g_1(x,y,z_{0},0,0) + \mathcal{O}(\epsilon,\delta), \\
\dot{Z} &= G(x,y,z_{0},0,0)+ Z H(x,y,z_{0},0,0) + \mathcal{O}(\epsilon,\delta).
\end{align}
\end{subequations}
In the singular limit of $\epsilon=0$, we obtain
\begin{subequations}\label{singular_proof}
\begin{align}
y &=F(x,z_{0})+ \mathcal{O}(\delta), \label{singular_proofa} \\
\dot {y} &= g_1(x,y,z_{0},0,0) + \mathcal{O}(\delta), \\
\dot{Z} &= G(x,y,z_{0},0,0)+ Z H(x,y,z_{0},0,0) + \mathcal{O}(\delta).
\end{align}
\end{subequations}
Differentiating \eqref{singular_proofa} with respect to $\tau$, we find
\begin{align*}
\dot{y}=F_{x}(x,z_{0})\dot{x}+ \mathcal{O}(\delta)
\end{align*}
which, together with the $(y,Z)$-subsystem of \eqref{singular_proof}, yields the projection of the reduced flow (in $\epsilon$) onto the $(x,Z)$-plane:
\begin{subequations}\label{red:slow_flow}
\begin{align}
\dot{x}  &=\frac{g_1(x,F(x,z_{0}),z_{0},0,0)}{F_{x}(x,z_{0})}+ \mathcal{O}(\delta), \\
\dot{Z}&= G(x,F(x,z_{0}),z_{0},0,0)+ Z H(x,F(x,z_{0}),z_{0},0,0) + \mathcal{O}(\delta).
\end{align}	
\end{subequations}
The important observation now is that Equation~\eqref{red:slow_flow} is partially decoupled for $\delta=0$. As a consequence, Assumption~\ref{A5} on the existence of a singular orbit $\Gamma_0$ actually implies that $g_1(x,F(x,z_{0}),z_{0},0,0)$ is non-zero. In other words, we may parametrize the reduced flow by the variable $x$. Hence, introducing $x$ as the independent variable in \eqref{red:slow_flow} and noting that $F_{x}(x,z_{0})\ne 0$ away from $L_i$, we obtain an ordinary differential equation
\begin{align}\label{ODE}
\frac{dZ}{dx} = p(x)Z + q(x) + O(\delta),
\end{align}
with
\begin{align*}
p(x) := \frac{H(x,F(x,z_{0}),z_{0},0,0)}{g_1(x,F(x,z_{0}),z_{0},0,0)} F_{x}(x,z_{0})
\ \text{and}\
q(x) := \frac{G(x,F(x,z_{0}),z_{0},0,0)}{g_1(x,F(x,z_{0}),z_{0},0,0)} F_{x}(x,z_{0}),
\end{align*}
which is linear with respect to $Z$ when $\delta=0$. 

In that limit, \eqref{ODE} can hence be solved exactly, with initial condition $Z(x_{\rm init})=Z_{\rm init}$:
\begin{equation}\label{ODE_soln}
Z(x,x_{\rm init},Z_{\rm init})=\big(e^{\int^{x}_{x_{\rm init}} p(s)ds  }\big) Z_{\rm init} +\left(\int^{x}_{x_{\rm init}} q(u) e^{ \int^{x}_{u} p(s) ds}du\right).
\end{equation}
Note that \eqref{ODE_soln} is precisely the affine map defined in Definition~\ref{assosc_def0} that is associated to the slow portion of $\Gamma_0^S$ or $\Gamma_0^L$ between $x_{\rm init}$ and $x$.
Now, we make use of \eqref{ODE_soln} to define a map that encodes the mixed-mode dynamics of our canonical system, Equation~\eqref{VFfast}. 
The discussion underneath Assumptions~\ref{A3} and \ref{A4} implies that the sought-after map will have two branches which describe oscillations with different amplitudes as we pass through $Z_{\rm init}=0$. Specifically, for $Z_{\rm init}<0$, we observe \textit{large-amplitude oscillations (LAOs)}, while for $Z_{\rm init}>0$, we have \textit{small-amplitude oscillations (SAOs)};
we hence proceed to define the following one dimensional piecewise affine map associated with \eqref{VFfast},
\begin{align} \label{PAM_proof}
M(Z_{\rm init})=
\begin{cases}
Z(x_{4},\hat{x}_{1},Z(x_{1},\hat{x}_{4},Z_{\rm init})) = M_{1}(Z_{init}) & \quad\text{if }Z_{\rm init}<0, \\
Z(x_{4},\hat{x}_{3},Z(x_{3},x_{2},Z_{\rm init}))=M_{2}(Z_{init}) & \quad\text{if }Z_{\rm init}>0,
\end{cases}								
\end{align}
where  $\hat{x}_{4}$, $\hat{x}_3$ and $\hat{x}_{1}$ are defined as in Figure~\ref{fig:critical_manifold}.
Given \eqref{ODE_soln}, we find the expressions for the affine maps defined in \eqref{PAM_proof} or, equivalently, in \eqref{PAM}:
\begin{multline}\label{PAM_LAO}
Z(x_{4}, \hat{x}_{1},Z(x_{1},\hat{x}_{4},Z_{\rm init})) = \Big({\rm e}^{\big(\int^{x_{1}}_{\hat{x}_{4}}+\int^{x_{4}}_{\hat{x}_{1}}\big)p(s)ds}\Big) Z_{\rm init} \\
+\bigg(\int^{x_{1}}_{\hat{x}_{4}} q(u){\rm e}^{\big(\int^{x_{1}}_{u}+\int^{x_{4}}_{\hat{x}_{1}}\big) p(s) ds}du + \int^{x_{4}}_{\hat{x}_{1}} q(u) {\rm e}^{\int^{x_{4}}_{u} p(s) ds}du\bigg)
\end{multline}
and
\begin{multline}\label{PAM_SAO}
Z(x_{4}, \hat{x}_{3},Z(x_{3},{x}_{2},Z_{\rm init})) = \Big({\rm e}^{\big(\int^{x_{3}}_{x_{2}}+\int^{x_{4}}_{\hat{x}_{3}}\big)p(s)ds}\Big) Z_{\rm init} \\
+\bigg(\int^{x_{3}}_{{x}_{2}} q(u){\rm e}^{\big(\int^{x_{3}}_{u}+\int^{x_{4}}_{\hat{x}_{3}}\big) p(s) ds}du + \int^{x_{4}}_{\hat{x}_{3}} q(u) {\rm e}^{\int^{x_{4}}_{u} p(s) ds}du\bigg).		
\end{multline}	

\begin{proposition}\label{return-limit-prop}
Let $\tilde{\Pi}_0$ be defined as in Section~\ref{statement of the results}. For each $Z\ne 0$, $\tilde{\Pi}_0$ is well-defined for $\delta>0$ sufficiently small; moreover, the limit $\tilde{\Pi}_{00}(Z)=\lim_{\delta\to 0}\tilde{\Pi}_0(Z)$ exists and is equal to the piecewise affine map given in \eqref{PAM_proof}, i.e., $\tilde{\Pi}_{00}(Z_{\rm init})=M(Z_{\rm init})$, with $Z_{\rm init}\ne 0$.
\end{proposition}
Proposition \ref{return-limit-prop} will be proved in Section~\ref{section-prop-proof}. In the proof, we will use an important observation made in Section~\ref{proof of first return map theorem}: the $\mathcal{Z}$-component of the return map $\Pi$ defined in Theorem \ref{firstreturnmap_thm} is a small $\epsilon$-perturbation of the return map induced by the reduced flow of \eqref{VFfast} near $\Gamma_0^L$ and $\Gamma_0^S$, respectively, provided we are below and above $B_{r}(P(L_2))$, respectively.

\section{Proof of main results}\label{section-proofs-all}

In this section, we present rigorous proofs for our main results, as introduced in Section~\ref{statement of the results}.

\subsection{Proof of Theorem~\ref{firstreturnmap_thm}}\label{proof of first return map theorem}

We first prove Theorem~\ref{firstreturnmap_thm}. To that end, we consider Equation~\eqref{VFfast} under Assumptions \ref{A1} through \ref{A5} to show that there exists an open neighborhood $\mathcal{U}$ of the intersection point $P_c$ of $P(L_2)$ and $P(L_4)$ such that, for all $\delta>0$ and $r>0$ small and fixed, the Poincar\'{e} map $\Pi: \mathcal{U}\setminus B_{r}(P(L_2))\subset\Delta\to\Delta$ induced by (\ref{VFfast}) is well-defined for $\epsilon$ sufficiently small. 

\begin{figure}[H] 
\includegraphics[width=12cm,center]{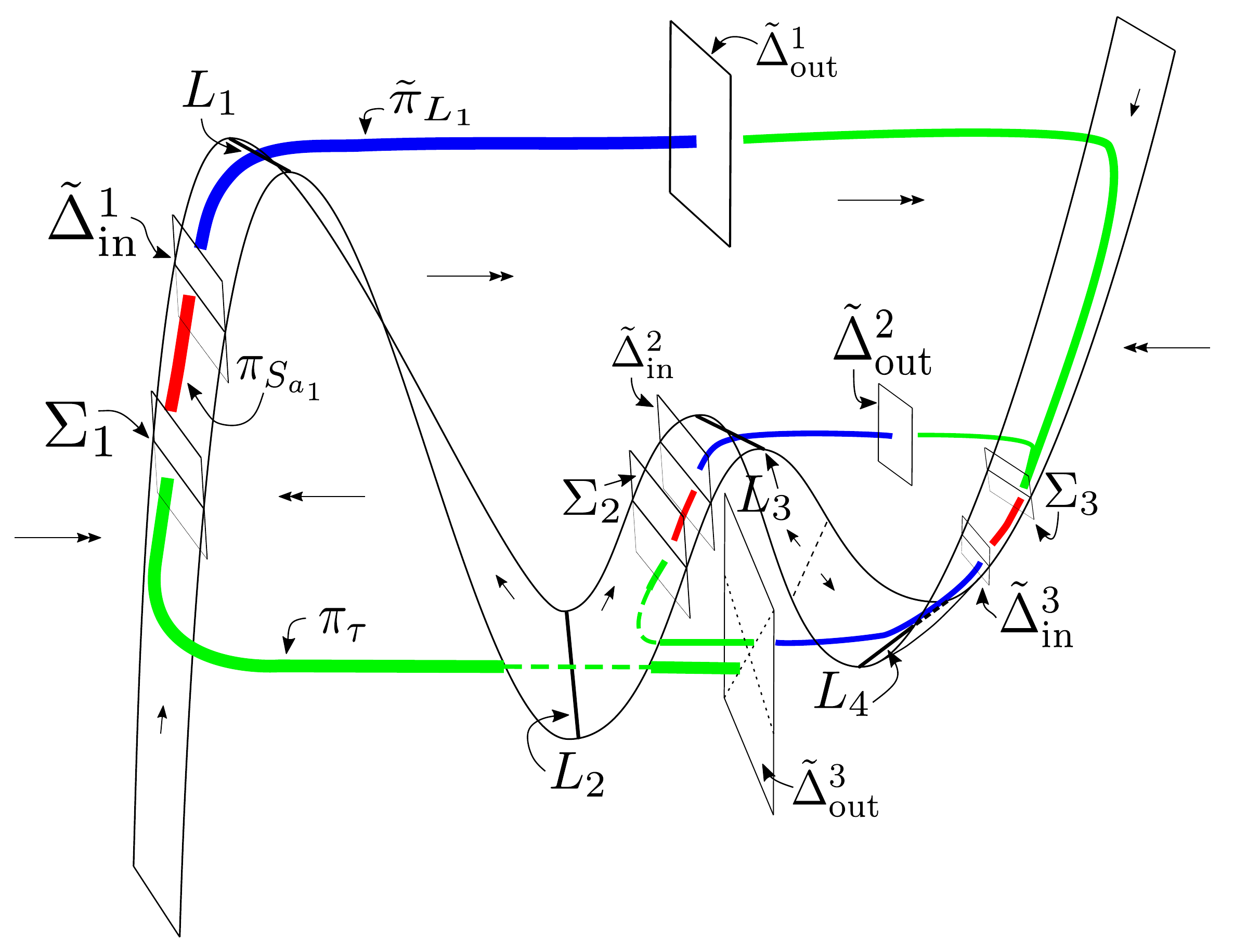}
\caption{The return map $\Pi$ induced by (\ref{VFfast}) on $\tilde{\Delta}^{3}_{out}$ is a composition of transition maps. The projection $P(L_2)$ divides $\tilde{\Delta}^{3}_{out}$ into two portions; in one portion, the fast flow is attracted to $\mathcal{S}_{a_1}$, whereas in the other, it tends to $\mathcal{S}_{a_2}$.}
\label{fig:Thm21-proof-sketch}
\end{figure}

Our proof is based on the techniques developed in \cite{SW04}, as indicated in Figure~\ref{fig:Thm21-proof-sketch}: ``fast" orbits of \eqref{VFfast} passing through $\Delta=\tilde{\Delta}^{3}_{\rm out}$ below the tubular neighbourhood $B_{r}(P(L_2))$ are attracted to $\mathcal{S}_{a_1}$, and therefore give rise to LAOs in the resulting mixed-mode time series; 
similarly, orbits passing through $\tilde{\Delta}^{3}_{\rm out}$ above $B_{r}(P(L_2))$ are attracted to $\mathcal{S}_{a_2}$, resulting in SAOs. Considered separately, each of these two cases can clearly be reduced to the return map studied in \cite{SW04}, for fixed $\delta>0$ and $r>0$. (Recall that we stay uniformly away from the fold line $L_2$.) We focus on the first case of LAOs here; the
second case, of SAOs, can be studied in an analogous fashion.

Fundamentally, we need to show that, for $r>0$ and $\delta>0$ sufficiently small, the flow of
\eqref{VFfast} stays close to the singular closed orbit $\Gamma_{0}^{L}$ such that the return 
map $\Pi$ exists for $\epsilon$ small.
Following \cite{SW04}, the map $\Pi$ is essentially composed of three different types of transition maps: $\pi_T$, $\pi_{S_{a_1}}$, and $\tilde{\pi}_{L_1}$, as illustrated in Figure~\ref{fig:Thm21-proof-sketch}. Here, the map $\pi_T$ is defined by following the fast flow towards the attracting portion $\mathcal{S}_{a_1}$ of $\mathcal{S}$, while $\pi_{S_{a_1}}$ describes the passage near $\mathcal{S}_{a_1}$ away from the fold line $L_1$; the study of $\pi_T$ and $\pi_{S_{a_1}}$ is based on Fenichel's standard GSPT. The map $\tilde{\pi}_{L_1}$, which describes the passage near the fold line $L_1$, is studied via geometric desingularisation, or ``blow-up". 
Let us now consider the ``half-return" map $\Pi_{H_\alpha}=\tilde{\pi}_{L_{1}} \circ \pi_{S_{a_{1}}} \circ \pi_{T}$ from the portion of $\mathcal{U}$ below $B_{r}(P(L_2))$ to a section $\tilde{\Delta}^{1}_{\rm out}$ transverse to $\Gamma_{0}^{L}$. Following Theorem 2 in \cite{SW04}, the half-return map $\Pi_{H_\alpha}$ is given by
\begin{equation*}
\Pi_{H_\alpha}(y,z)=(\mathcal{Y}_\alpha(y,z,\epsilon,\delta),\mathcal{Z}_\alpha(y,z,\epsilon,\delta)),
\end{equation*}
with $\mathcal{Z}_\alpha(y,z,\epsilon,\delta)=\mathcal{Z}_\alpha(y,z,\delta)+O(\epsilon\ln\epsilon)$ where $\mathcal{Z}_\alpha(y,z,\delta)$ is defined by following the orbit of the reduced flow on the attracting portion $\mathcal{S}_{a_1}$ between the $\omega$-limit of the point $(y,z)$ and $L_1$. Moreover, we have $\mathcal{Y}_\alpha(y,z,\epsilon,\delta)=\phi_\alpha(\mathcal{Z}_\alpha(y,z,\delta),\delta)+o(1)$, where $y=\phi_\alpha(z,\delta)$ describes the projection of the fold $L_1$ onto $\tilde{\Delta}^{1}_{\rm out}$ and where the $o(1)$-term tends uniformly to zero as $\epsilon\to 0$.

The half-return map $\Pi_{H_\beta}=(\mathcal{Y}_\beta,\mathcal{Z}_\beta)$ from $\tilde{\Delta}^{1}_{\rm out}$ back to $\Delta$ can be studied in a similar fashion, as a composition of transition maps that are of the same type as in $\Pi_{H_\alpha}$. Combining the two, the return map $\Pi=\Pi_{H_\beta}\circ\Pi_{H_\alpha}$, which is defined in the region below $B_{r}(P(L_2))$, can be written as 
\begin{equation*}
\Pi(y,z)=(\mathcal{Y}(y,z,\epsilon,\delta),\mathcal{Z}(y,z,\epsilon,\delta)),
\end{equation*} 
with $\mathcal{Z}(y,z,\epsilon,\delta)=\mathcal{Z}(y,z,\delta)+O(\epsilon\ln\epsilon)$, where 
\begin{equation}\label{return-reduced}
\mathcal{Z}(y,z,\delta)=\mathcal{Z}_\beta(\phi_\alpha(\mathcal{Z}_\alpha(y,z,\delta),\delta),\mathcal{Z}_\alpha(y,z,\delta),\delta)
\end{equation}
is the return map defined by the reduced flow on $\mathcal{S}_{a_1}$ and $\mathcal{S}_{a_3}$. We can also conclude that the function $\mathcal{Y}$ has the property given in Theorem \ref{firstreturnmap_thm}, which completes the proof. 

\subsection{Proof of Proposition~\ref{return-limit-prop}}\label{section-prop-proof} 

Next, we prove Proposition~\ref{return-limit-prop}. Recall that the map
\[
\tilde{\Pi}_0(Z_{\rm init})= \frac{\mathcal{Z}_{1}(z_0+\delta Z_{\rm init},\delta) -z_0}{\delta},
\]
with $\mathcal{Z}_{1}(z,\delta) = \mathcal{Z}(\phi_{L_4}(z,\delta),z,0,\delta)$, is defined for $z\ne z_0(\delta)$; cf.~Section~\ref{statement of the results}. If $Z_{\rm init}\ne 0$ is fixed, then $z=z_0+\delta Z_{\rm init}\ne z_0(\delta)$ for $\delta>0$ sufficiently small due to $z_{0}'(0)=0$; see Assumption~\ref{A3}. Thus, $\tilde{\Pi}_0$ is well-defined for $Z_{\rm init}\ne 0$ provided that $\delta>0$ is small.

First, let us consider $Z_{\rm init}<0$ and fixed. Then, we have that $z=z_0+\delta Z_{\rm init}<z_0(\delta)$ for $\delta>0$ small, i.e., the point $(\phi_{L_4}(z,\delta),z)$ is attracted to $\mathcal{S}_{a_1}$; see again Assumption~\ref{A3}. We therefore observe LAOs and $\mathcal{Z}_{1}(z,\delta) = \mathcal{Z}(\phi_{L_4}(z,\delta),z,\delta)$, where the function $\mathcal{Z}(y,z,\delta)$ is defined in (\ref{return-reduced}). 
Now, we note that the system in (\ref{red:slow_flow}) is obtained by applying the coordinate transformation in (\ref{eq:rescale}) to the reduced flow in \eqref{VFreduced} in $(x,z)$-space, where $g_2$ is given in (\ref{eq:g2}). It follows that the orbit of (\ref{red:slow_flow}) which is initiated at $(\hat{x}_{4}(z_0+\delta Z_{\rm init},\delta), Z_{\rm init})$, with $F(\hat{x}_{4}(z,\delta),z,0,\delta)=\phi_{L_4}(z,\delta)$ -- see Figure \ref{fig:critical_manifold} -- intersects the projection of the fold line $L_1$ onto the $(x,Z)$-space in $(x_{1}(z_0+\delta Z_\alpha,\delta),Z_\alpha)$, where 
\begin{equation*}
Z_\alpha=\frac{\mathcal{Z}_\alpha(\phi_{L_4}(z_0+\delta Z_{\rm init},\delta),z_0+\delta Z_{\rm init},\delta)-z_0}{\delta}.
\end{equation*}
(Here, $\mathcal{Z}_\alpha$ is defined as in Section \ref{proof of first return map theorem}.)
Thus, $Z_\alpha$ converges to $Z(x_{1},\hat{x}_{4},Z_{\rm init})$ as $\delta\to 0$, with $Z(x,x_{\rm init},Z_{\rm init})$ given in (\ref{ODE_soln}), where we denote by $\hat{x}_{4}$ and $x_{1}$, respectively, the limit of $\hat{x}_{4}(z,\delta)$ and $x_{1}(z,\delta)$, respectively, as $(z,\delta)\to(z_0,0)$. Here, we have used the fact that $\Gamma_0^L$ is located in the plane $\{z=z_0\}$ with $\delta=0$ and, thus, that $\delta Z_\alpha\to 0$ as $\delta\to 0$ in $x_1$. Moreover, we have exploited our observation in Section~\ref{sec:derivePAM} that (\ref{red:slow_flow}) is a $\delta$-perturbation of a linear (in $Z$) differential equation.

Similarly, the orbit of (\ref{red:slow_flow}) that is initiated at $(\hat{x}_{1}(z_0+\delta Z_\alpha,\delta),Z_\alpha)$, with $F(\hat{x}_{1}(z,\delta),z,0,\delta)=\phi_{\alpha}(z,\delta)=F({x}_{1}(z,\delta),z,0,\delta)$, again by Figure \ref{fig:critical_manifold}, intersects the $(x,Z)$-projection of the fold line $L_4$ in $(x_{4}(z_0+\delta Z_\beta,\delta),Z_\beta)$, with 
\begin{equation*}
Z_\beta=\frac{\mathcal{Z}_\beta(\phi_{\alpha}(z_0+\delta Z_{\alpha},\delta),z_0+\delta Z_{\alpha},\delta)-z_0}{\delta}.
\end{equation*}
We therefore conclude that $Z_\beta$ converges to (\ref{PAM_LAO}) as $\delta\to 0$. (As above, 
we use that $\delta Z_\beta\to 0$ as $\delta\to 0$ in $x_4$.) Now, it suffices to note that $\tilde{\Pi}_0(Z_{\rm init})=Z_\beta$, from (\ref{return-reduced}).

The case where $Z_{\rm init}>0$ can be studied in a similar fashion to show that $\tilde{\Pi}_0(Z_{\rm init})$ tends to (\ref{PAM_SAO}) as $\delta\to 0$, as claimed, which completes the proof.

\subsection{Proof of Theorem~\ref{association_thm}}\label{proof of assosciation theorem}

For the sake of simplicity and readability, we first prove Theorem~\ref{association_thm} for a MMO with signature $1^0$; then, we will indicate how the proof can be extended to the general case, i.e.,
to MMOs with signature $L_1^{s_1}L_2^{s_2}\cdots L_k^{s_k}$, with $k\ge 1$ integer. 

Thus, we suppose that the PAM in \eqref{PAM} which is associated with the vector field in \eqref{VFfast} has a stable periodic orbit that undergoes one LAO, i.e., that $M(Z^\ast)=Z^\ast$ for $Z^\ast<0$ and $a_{11}<1$, where $a_{11}$ is the coefficient of $Z$ in $M$. Our goal is to prove that \eqref{VFfast} has a stable periodic orbit with one LAO for $\epsilon,\delta>0$ small. Clearly, periodic orbits for \eqref{VFfast} correspond to fixed points of the first return map $\Pi$ defined in Theorem~\ref{firstreturnmap_thm}. It can easily be seen that $(y,z)$ is a solution of $\Pi(y,z)-(y,z)=(0,0)$ if and only if $(y,Z)$, with $z=z_0+\delta Z$, is a solution of 
\begin{align}\label{fixed points-1}
\begin{split}
\phi_{L_4}(z_0+\delta\tilde{Z}(y,Z,0,\delta),\delta)+\mathcal{E}(y,z_0+\delta Z,\epsilon,\delta)-y &=0, \\
\tilde{Z}(y,Z,\epsilon,\delta)-Z &=0,
\end{split}
\end{align}
where 
\[
\tilde{Z}(y,Z,\epsilon,\delta):=\frac{\mathcal{Z}(y,z_0+\delta Z,\epsilon,\delta)-z_0}\delta
\]
and $\phi_{L_4}$, $\mathcal{Z}$, and $\mathcal{E}$ are defined as in Theorem \ref{firstreturnmap_thm}. Using the Implicit Function Theorem, we show that the system in \eqref{fixed points-1} has a unique solution $(y^\ast_{\epsilon,\delta},Z^\ast_{\epsilon,\delta})$ for $\epsilon>0$ and $\delta>0$ sufficiently small, with $(y^\ast_{\epsilon,\delta},Z^\ast_{\epsilon,\delta})$ close to $(\phi_{L_4}(z_0,0),Z^\ast)$. (An alternative approach is outlined in Remark~\ref{approach-2}.) Note that $\tilde{Z}(\phi_{L_4}(z_0+\delta Z,\delta),Z,0,\delta)=\tilde{\Pi}_0(Z)$, where $\tilde{\Pi}_0(Z)$ tends to $M(Z)$ as $\delta\to 0$, by Proposition~\ref{return-limit-prop}. More generally, we have $\tilde{Z}(y,Z,0,\delta)\to M(y,Z)$ as $\delta\to 0$, where $M(y,Z)$ is a PAM as in \eqref{PAM_proof} or, equivalently, in \eqref{PAM}, with $\hat{x}_4$ and $x_2$ depending on $y$. (This follows easily from the proof of Proposition \ref{return-limit-prop}.) Now, letting $\epsilon \to 0$ and then $\delta\to 0$, the system in \eqref{fixed points-1} reduces to
\begin{align}\label{fixed points-2}
\begin{split}
\phi_{L_4}(z_0,0)-y &=0, \\
M(y,Z)-Z &=0.
\end{split}
\end{align}
Since $Z=Z^*$ is a fixed point of $M(Z)$ -- or, equivalently, of $M(\phi_{L_4}(z_0,0),Z)$, it follows that $(y,Z)=(\phi_{L_4}(z_0,0),Z^*)$ is a solution of \eqref{fixed points-2}. The Jacobian determinant of the left-hand side in \eqref{fixed points-2} evaluated at this solution is $1-a_{11}\ne 0$, where we note that $M'(Z^\ast)=a_{11}$ because $Z^\ast<0$. The Implicit Function Theorem now implies the existence of a solution $(y^\ast_{\epsilon,\delta},Z^\ast_{\epsilon,\delta})$ of \eqref{fixed points-1} for $\epsilon,\delta>0$ small. Thus, $(y,z)=(y^\ast_{\epsilon,\delta},z_0+\delta Z^\ast_{\epsilon,\delta})$ is a fixed point of $\Pi$. It is clear that the corresponding periodic orbit is stable. This completes the proof.

In the general case, where the given MMO has signature $L_1^{s_1}L_2^{s_2}\cdots L_k^{s_k}$, 
 we have to study fixed points of the $\kappa$-th iterate of the first return map $\Pi$, where $\kappa:=\sum^{k}_{i=1}{\big( L_{i}+s_{i}}\big)$.
In the limit of $\epsilon=0=\delta$, the $\kappa$-th iterate of $\Pi$ can be written as $(y,Z)\to (\phi_{L_4}(z_0,0),M^{\kappa-1}(M(y,Z)))$ in $(y,Z)$-coordinates. 
The Jacobian determinant of the corresponding system $\{\phi_{L_4}(z_0,0)-y=0,M^{\kappa-1}(M(y,Z))-Z=0\}$ is then equal to $1-a_{11}^L a_{21}^s$, with $L=\sum^{k}_{i=1} L_{i}$ and $s=\sum^{k}_{i=1} s_{i}$. 
Since we supposed that $M^\kappa(Z^\ast)=Z^\ast$ for some $Z^\ast<0$ with $a_{11}^L a_{21}^s<1$ (stability), the result easily follows.

\begin{remark}\label{approach-2}
Alternatively, Theorem~\ref{association_thm} can be proved via the approach taken in \cite{SW04}. 
For $\delta>0$ small, the first return map $\Pi$ from Theorem~\ref{firstreturnmap_thm} contracts its domain to the curve $y=\phi_{L_4}(z,\delta)$, in the limit as $\epsilon\to 0$. 
Following Theorem 3 in \cite{SW04}, $\Pi$ admits a one-dimensional attracting invariant manifold $y=m_{\epsilon,\delta}(z)$; the dynamics of $\Pi$ on that manifold is given by $\Pi_0$ in the limit of $\epsilon\to 0$, with $\Pi_0$ as defined underneath Theorem~\ref{firstreturnmap_thm}. 
In $(y,Z)$-coordinates, $\Pi_0$ is given by the PAM $M(Z)$ for $\delta\to 0$; see 
Proposition~\ref{return-limit-prop}. 
Now, it suffices to note that hyperbolic fixed points of $M$ persist under perturbation 
of $M$ in $\delta$ -- which gives $\Pi_0$ -- and, subsequently, under perturbation of $\Pi_0$ in $\epsilon$. Thus, we find a fixed point of the one-dimensional map $\Pi_{m_{\epsilon,\delta}}$; 
the $\kappa$-th iterate of $\Pi$ can be studied in a similar fashion.
\end{remark}

\subsection{Proof of Theorem~\ref{existence_thm}}\label{Proof of existence theorem}

To prove Theorem~\ref{existence_thm}, we introduce a specific sub-family of slow-fast systems
of the form in \eqref{VFfast} that satisfies Assumptions~\ref{A1} through \ref{A5}. Then, we will
show that a given PAM $M$ can be associated to a representative system from that family. 
Specifically, we take
\[
F(x,z,\epsilon,\delta)=F(x,z),\ g_1(x,y,z,\epsilon,\delta)=J(x),\ \text{and}\ g_2(x,y,z,\epsilon,\delta)=\delta G(x)+z H(x)
\]
in \eqref{VFfast}. In particular, we take $F(x,z)$ to be a polynomial of degree $9$ in $x$,
restricted to $(x,z)\in (-3,2)\times (-1,1)$; moreover, we choose the functions $G(x)$, $H(x)$, $J(x)$, 
and $Q(x)$ such that the integrals to be evaluated in \eqref{PAM_LAO} and \eqref{PAM_SAO} are as
simple as possible, with convenient substitutions inside the integrands. Also, for simplicity, we take
$z_0=0$.

In sum, we hence have
\begin{subequations}\label{VF_represent}
\begin{align}
\begin{split}
F(x,z) &=a_{9}x^{9}+\sum_{k=2}^8a_{k}(z)x^{k},\ \text{with} \ a_{9}= \tfrac{184180}{67741437},
\ a_{8}(z)=\tfrac{1}{8}\big(\tfrac{138135}{90321916}z+\tfrac{3558512}{22580479}\big), \\
a_{7}(z) &=\tfrac{1}{7}\big(\tfrac{751493}{90321916}z+\tfrac{212863}{22580479}\big), \  a_{6}(z)=-\tfrac{1}{6}\big(\tfrac{2793109}{361287664}z+\tfrac{23361467}{22580479}\big), \\
a_{5}(z) &=-\tfrac{1}{5}\big(\tfrac{10284179}{180643832}z+\tfrac{1224990}{22580479}\big),\
a_{4}(z)=\tfrac{1}{4}\big(\tfrac{2417921}{45160958}z+\tfrac{64963913}{22580479}\big), \\
a_{3}(z) &=\tfrac{1}{3}\big(\tfrac{45620545}{361287664}z+\tfrac{459587}{22580479}\big), \ \text{and} \
a_{2}(z)=-\tfrac{1}{2}\big(\tfrac{1}{8}z+2\big),
\end{split}\\
J(x)&=\frac{1}{2}-x, \ \rho(x)=p+x+qx^2, \ Q(x)=\int^{x}_{0} \rho(s) \frac{\partial F}{\partial s}(s, z_0)ds,  \\
G(x)&=\Big[\kappa+\lambda\Big(\frac{\alpha Q(x)^2}{2}+\beta Q(x)\Big)\Big](\alpha Q(x)+\beta)\rho(x)J(x), \ \text{and} \\
H(x)&=\rho(x)\bigg[\int^{x}_{0}\Big(\rho(s) \frac{\partial F}{\partial s}(s, z_0)ds\Big)\alpha+ \beta\bigg]J(x).
\end{align}
\end{subequations}

While the choices in \eqref{VF_represent} seem far from simple at first glance, they are made for
the sole purpose of simplifying the requisite calculations that follow. (A related system will also underlie the numerical simulations presented in the next Section~\ref{sec:numerics}; although that system will mostly be identical to the one in \eqref{VF_represent}, the definition of the function $\rho(x)$ will differ for computational efficiency.)

By Definition~\ref{PAM}, the given PAM $M(Z)$ is determined by the coefficients $a_{ij}$, for $i,j=1,2$. We will prove that there exists a fast-slow system of the specific form
in \eqref{VF_represent} which is associated with $M$; to that end, we need to show that the system of equations
\begin{align}\label{eq:aij}
a_{ij}=f_{ij}(\alpha,\beta,\kappa,\lambda,p,q),\quad\text{with }i,j=1,2,
\end{align}
has at least one solution $(\alpha_\ast,\beta_\ast,\kappa_\ast,\lambda_\ast,p_\ast,q_\ast)$ which fully determines the vector field in \eqref{VF_represent}.
(Here, the notation $f_{ij}$ is shorthand for
the right-hand sides in the definition of $a_{ij}$ in Definition~\ref{PAM}.)

In a first step, we note that $a_{11}$ and $a_{21}$ depend on $(\alpha,\beta,p,q)$ only, i.e., that
\begin{align*}
a_{11} &=g_{11}(\alpha,\beta,p,q)\quad\text{and} \\
a_{21} &=g_{21}(\alpha,\beta,p,q)
\end{align*}
for some new functions $g_{11}$ and $g_{21}$, as well as that $a_{11}$ and $a_{21}$ are positive by definition:
\begin{align*}
a_{11} &= g_{11}(\alpha,\beta,p,q)= \exp\big(\mathcal{A}_{11}(p,q)\alpha + \mathcal{A}_{12}(p,q)\beta\big)\quad\text{and} \\
a_{21} &= g_{21}(\alpha,\beta,p,q)= \exp\big(\mathcal{A}_{21}(p,q)\alpha + \mathcal{A}_{22}(p,q)\beta\big)
\end{align*}
for some functions $\mathcal{A}_{ij}(p,q)$ which are, in fact, polynomial in $p$ and $q$. 
Taking logarithms, we find a linear system in the unknowns $(\alpha,\beta)$ whose principal matrix has determinant $\Delta(p,q)$. 
With the aid of the computer algebra package {\sc Maple}, we compute $\Delta$ to be a polynomial 
of degree $3$ with positive coefficients. Restricting to the parameter domain $\{p>0,q>0\}$, we can hence safely assume that $\Delta$ is non-zero and, hence, that the above system has a solution
\[
\alpha = \alpha_\ast(p,q,a_{11},a_{21})\quad\text{and}\quad\beta = \beta_\ast(p,q,a_{11},a_{21}).
\]

The next part of the proof is more intricate, and again relies on symbolic computation in {\sc Maple}.  Substituting the above expressions for $(\alpha,\beta)=(\alpha_\ast,\beta_\ast)$ into 
\eqref{eq:aij}, we obtain
\begin{align}\label{eq:aij2}
a_{12} &=f_{12}(\alpha_\ast,\beta_\ast,\kappa,\lambda,p,q) \\
a_{22} &=f_{22}(\alpha_\ast,\beta_\ast,\kappa,\lambda,p,q),
\end{align}
which is a linear system in $(\kappa,\lambda)$ whose principal matrix has determinant $\tilde\Delta(p,q,a_{11},a_{21})$. The expression for $\tilde\Delta$ can be written as
\begin{align*}
\tilde\Delta(p,q,a_{11},a_{21})=\frac{N( p,q,a_{11},a_{21},\ln a_{11},\ln a_{21},{\rm e}^{E_1},{\rm e}^{E_2})}{\Delta( p,q)},
\end{align*}
where $\Delta$ is defined as above, the exponents of the exponential terms ${\rm e}^{E_i}$ are of the form $E_i=E_i(\alpha,\beta,p,q)$, and $N$ is polynomial in all its $8$ arguments. (In fact, 
$N$ has degree 1 with respect to $e^{E_1}$ and $e^{E_2}$.)
It now suffices to show that, for each choice of $(a_{11},a_{21})$, there is at least one choice of $(p,q)$, with $p>0$ and $q>0$, for which $N$ is non-zero. Given the complex algebraic form of $N$, that is a cumbersome task. However, it suffices to argue that almost any choice of $(p,q)$ will be 
admissible.
\medskip

We will outline that argument here. First, we write $N = N_0 + N_1 e^{E_1} + N_2 e^{E_2}$, where 
each $N_i$ is a polynomial expression in $(p,q,a_{11},a_{21})$. Using {\sc Maple}, we verify that
\[
\lim_{p\to\pm \infty}\frac{E_i}{p} = \frac{R_i(q,\ln a_{11},\ln a_{12})}{S(q)},
\]
for some strictly positive degree-$2$ polynomial $S$ and some degree-$1$ polynomials (in $q$) $R_1$ and $R_2$. As there is only one choice for $q$ where the asymptotics of $E_1$ coincides with that of $E_2$, we restrict to the generic case where the two limits are strictly different. Then, there are six possibilities,
\begin{subequations}
\begin{align}
\textstyle\lim_{p\to\infty}\frac{E_1}{p}	&< \textstyle\lim_{p\to  \infty}\frac{E_2}{p} < 0, \\
\textstyle\lim_{p\to\infty}\frac{E_1}{p}	&<0 \leq\textstyle\lim_{p\to\infty}\frac{E_2}{p},\quad\text{and} \\
0 &\leq\textstyle\lim_{p\to\infty}\frac{E_1}{p}<\textstyle\lim_{p\to\infty}\frac{E_2}{p},
\end{align}
\end{subequations}
as well as the three possibilities obtained by swapping $E_1$ and $E_2$. Let us consider the third case as an example: in that scenario, as $p\to\infty$, the contributions of $N_0$ and $N_1$ in $N$ become negligible, and it suffices to see whether or not one can find $q$ for which $N_2$ is non-zero.  Equally, in the first two scenarios, we see that the contribution of $N_1$ becomes significant in the limit as $p\to-\infty$.
It now suffices to observe that both $N_1/p$ and $N_2/p$ are asymptotic to a quadratic polynomial in $q$ for large $|p|$ and, hence, that there are many choices of $(p,q)$ for which these expressions are non-zero. Hence, at least for $|p|$ sufficiently large, one can solve Equation~\eqref{eq:aij2} for
$(\kappa,\lambda)$, which, in sum, gives a solution $(\alpha_\ast,\beta_\ast,\kappa_\ast,\lambda_\ast)$ to \eqref{eq:aij}.  Hence, generically, given a PAM $M$, one can choose
$(p_\ast,q_\ast)$ such that there exists a slow-fast vector field within the family defined by \eqref{VF_represent} to which $M$ is associated. This completes the proof.
\begin{remark}
In practice, one would not take $|p|$ too large, as that would introduce another layer of time scale separation in the system.
\end{remark}

\section{Numerical verification}\label{sec:numerics}

Finally, in this section, we give a numerical verification of two of our main results, Theorems~\ref{association_thm} and \ref{existence_thm}.
To that end, we consider the family of one-dimensional PAMs of the form 
\begin{equation} \label{eq:PAM_prototype}
M(Z)=\begin{cases}
a_{11} Z + a_{12} & \quad\text{for }Z<0, \\
a_{21} Z + a_{22} & \quad\text{for }Z>0,
\end{cases}	
\end{equation} 
where $a_{ij}=f_{ij}(\alpha,\beta,\kappa,\lambda,p,q)$ with $i,j=1,2$, as introduced in~\eqref{eq:aij}.
For the calculations of the integrals appearing in (\ref{PAM_LAO}) and (\ref{PAM_SAO}), we require the following $x$-values, which are obtained from (\ref{VF_represent}) with $z_0=0$:
\[
\hat{x}_{4}= -\tfrac{5}{2},\ x_{1}=-2,\ x_{2}=-1,\ x_{3}=0,\ x_{4}=1,\ \hat{x}_{3}=\frac{3}{2},\ \text{and}\ \hat{x}_{1}=\frac{8}{5};
\]
see Figure \ref{fig:critical_manifold}. 
Next, and as outlined in the proof of Theorem~\ref{existence_thm} in Section~\ref{Proof of existence theorem}, we have to choose a suitable function $\rho$ in \eqref{VF_represent}. Rather than taking
$\rho$ within the family specified there, we pick the numerically more convenient function
\begin{equation}\label{eq:rho}
\rho(x) = \bigg( \frac{552540}{22580479}x^4+\frac{2453432}{22580479}x^3-\frac{4141461}{22580479}x^2-\frac{11520033}{22580479}x+1 \bigg)^{-1}.
\end{equation}
The choice in \eqref{eq:rho} allows us to determine the four pivotal quantities
$\alpha$, $\beta$, $\kappa$, and $\lambda$ in \eqref{VF_represent}, in agreement with our
expectation that a wide range of functions $\rho(x)$ will yield an admissible solution 
$(\alpha_\ast,\beta_\ast,\kappa_\ast,\lambda_\ast)$. That solution then specifies a
three-dimensional slow-fast system from the family determined by \eqref{VF_represent} 
that is associated to the given PAM in \eqref{eq:PAM_prototype}.

Below, we showcase a number of examples which verify that the resulting mixed-mode time series
in that system have identical signature to the corresponding periodic orbits for the PAM $M$, thus 
verifying Theorem~\ref{association_thm}. Here, we note that the functions $\rho$ and $J$ in 
\eqref{VF_represent} are independent of $\alpha$, $\beta$, $\kappa$ and $\lambda$, and that they 
hence do not change with the signature. The functions $Q$, $G$, and $H$, on the other hand,
are signature-dependent.

\begin{remark}
Given that $F(x,z)$ in \eqref{VF_represent} is a ninth-degree polynomial in $x$, 
we rescaled $x$ and $y$ as
\[
x\mapsto\frac27x \quad \text{and} \quad y\mapsto\frac32y
\]
in our visualisation in order to restrict the area of interest in $x$ to the interval $[-1,1]$. 
(We note that $z$ remains rescaled to $Z$, as defined by the change of coordinates in \eqref{eq:rescale} which was used in both Sections~\ref{statement of the results} and \ref{sec:derivePAM}).
\end{remark}

\subsection{Examples: MMOs of various signatures}

In a first step, we fix the coefficients $a_{11}$, $a_{21}$, and $a_{22}$ in the definition of the
PAM in \eqref{eq:PAM_prototype}, varying only $a_{12}$ as the ``bifurcation parameter". In Table~\ref{table:signatures_PAMS_abkl} below, we list two sequences of mixed-mode signatures that are obtained upon variation of $a_{12}$, with $a_{11}$, $a_{21}$, and $a_{22}$ fixed as stated there. For completeness, and to illustrate the two-way correspondence established in Theorems~\ref{association_thm} and \ref{existence_thm}, we also give the corresponding pivotal quantities $\alpha$, $\beta$, $\kappa$ and $\lambda$ in the definition of the associated vector field in \eqref{VF_represent}.
(We note that, given $a_{11}$, $a_{21}$, and $a_{22}$, $\alpha$ and $\beta$ do not change as $a_{12}$ is varied, in contrast to $\kappa$ and $\lambda$, as is to be expected from the proof of Theorem~\ref{existence_thm}.)

In particular, we thus observe an unfolding of a ``regular" sequence of signatures which are either 
of the form $\{1^s\}$ or $\{L^1\}$ in the bifurcation parameter $a_{12}$.
A selection of (periodic) MMO orbits, both for the PAM in \eqref{eq:PAM_prototype} and the associated vector field, is illustrated graphically in the figures below. 	We emphasise that we observe the same signature in all three (state) variables $x$, $y$, and $z$ in \eqref{VFfast}, which is due to
the geometry of the underlying critical manifold $\mathcal{S}$; see Figure~\ref{fig:MMO1-3} 
below, where we highlight the signature $1^3$ as one particular example.

\begin{table}[H]
	\setlength{\tabcolsep}{2pt}
		\centering 
		\begin{tabular}{|c|c|c|c|c|c|c|c|c|}
			\hline
			\hline                        
\tiny{\textbf{Signature}}&$a_{11}$&$a_{12}$&$a_{21}$&$a_{22}$	&\textbf{$\alpha$}&\textbf{$\beta$}&\textbf{$\kappa$}	&\textbf{$\lambda$}	\\
			\hline
			\hline
			
	 	$1^1$	&0.3		&1		&0.9		&-2		&0.8743		&0.0240		&27.	2674		&-64.5764			\\[2pt]
\hline  $1^2$	&0.3		&3		&0.9		&-2		&0.8743		&0.0240		&28.2364		&-73.1866			\\[2pt]
\hline	$1^3$	&0.3		&7		&0.9		&-2		&0.8743		&0.0240		&30.1744		&-90.4070			\\[2pt]
\hline	$1^4$	&0.3		&10		&0.9		&-2		&0.8743		&0.0240		&31.6279		&-103.3223			\\[2pt]
\hline	$1^5$	&0.3		&12		&0.9		&-2		&0.8743		&0.0240		&32.5969		&-111.9325			\\[2pt]
\hline	$1^6$	&0.3		&15		&0.9		&-2		&0.8743		&0.0240		&34.0504		&-124.8478			\\[2pt]
\hline	$1^7$	&0.3		&20		&0.9		&-2		&0.8743		&0.0240		&36.4729		&-146.3733			\\[2pt]
\hline	$1^8$	&0.3		&25		&0.9		&-2		&0.8743		&0.0240		&38.8954		&-167.8987			\\[2pt]
\hline
\hline
\hline	$1^1$	&0.9		&3		&0.4		&-3		&-0.5065		&1.0238		&3.2091		&-7,7202			\\[2pt]
\hline	$2^1$	&0.9		&1.5		&0.4		&-3		&-0.5065		&1.0238		&3.9766		&-4,4118			\\[2pt]
\hline	$3^1$	&0.9		&1		&0.4		&-3		&-0.5065		&1.0238		&4.2325		&-3.3088			\\[2pt]
\hline	$4^1$	&0.9		&0.7		&0.4		&-3		&-0.5065		&1.0238		&4.3860		&-2.6471			\\[2pt]
\hline	$5^1$	&0.9		&0.5		&0.4		&-3		&-0.5065		&1.0238		&4.4883		&-2.2059			\\[2pt]
\hline	$6^1$	&0.9		&0.4		&0.4		&-3		&-0.5065		&1.0238		&4.5395		&-1.9853			\\[2pt]
\hline	$7^1$	&0.9		&0.3		&0.4		&-3		&-0.5065		&1.0238		&4.6162		&-1.7647			\\[2pt]
\hline	$8^1$	&0.9		&0.25	&0.4		&-3		&-0.5065		&1.0238		&4.6418		&-1.6544			\\[2pt]
\hline			
    \end{tabular}
    \caption{Signatures of the form $L^1$ and $1^s$ generated by (\ref{VF_represent}) and the associated PAM in \eqref{eq:PAM_prototype}.}
    \label{table:signatures_PAMS_abkl}
	\end{table}	
			
\begin{figure}[H]
\includegraphics[width=12cm,center]{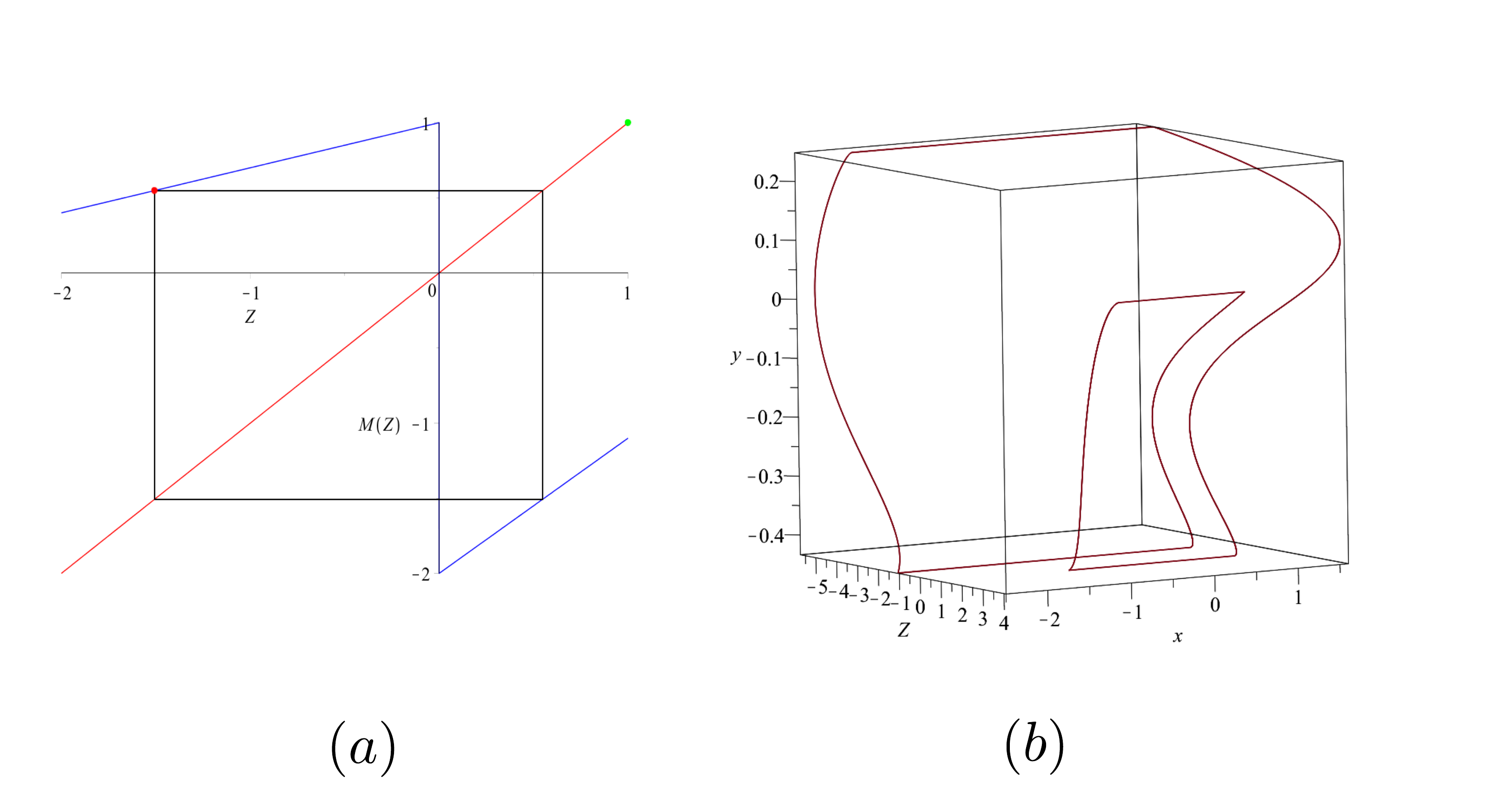} 
\caption{MMO of signature $1^1$: (a) piecewise affine map and (b) associated three-dimensional slow-fast system.}
\label{fig:MMO1-1}
\end{figure}

\begin{figure}[H]
\includegraphics[width=15cm,center]{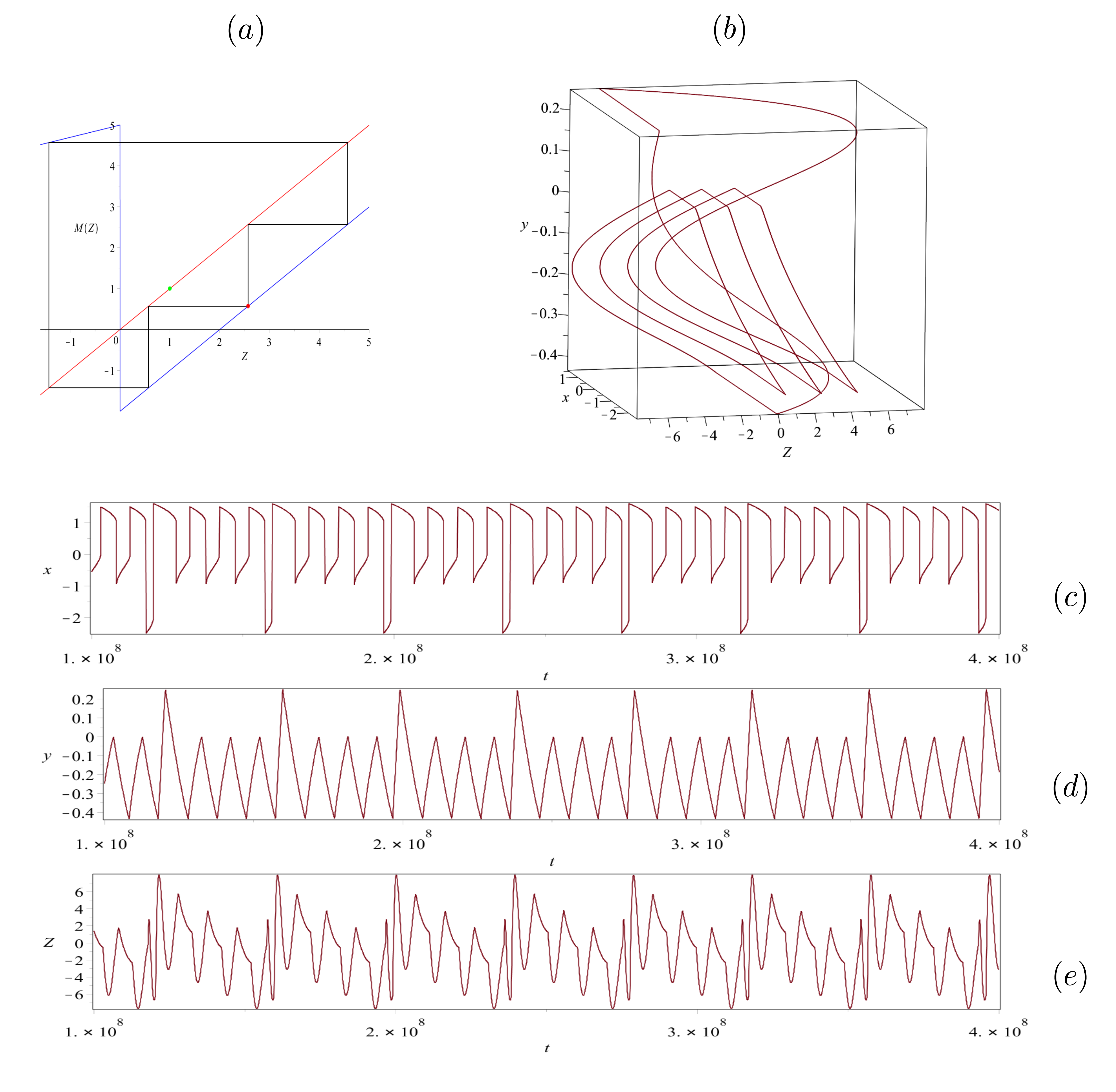} 
\caption{MMO of signature $1^3$: (a) piecewise affine map and (b) associated three-dimensional slow-fast system, with time series of $x$, $y$, and $Z$ in (c) through (e).}
\label{fig:MMO1-3}
\end{figure}

\begin{figure}[H]
\includegraphics[width=12cm,center]{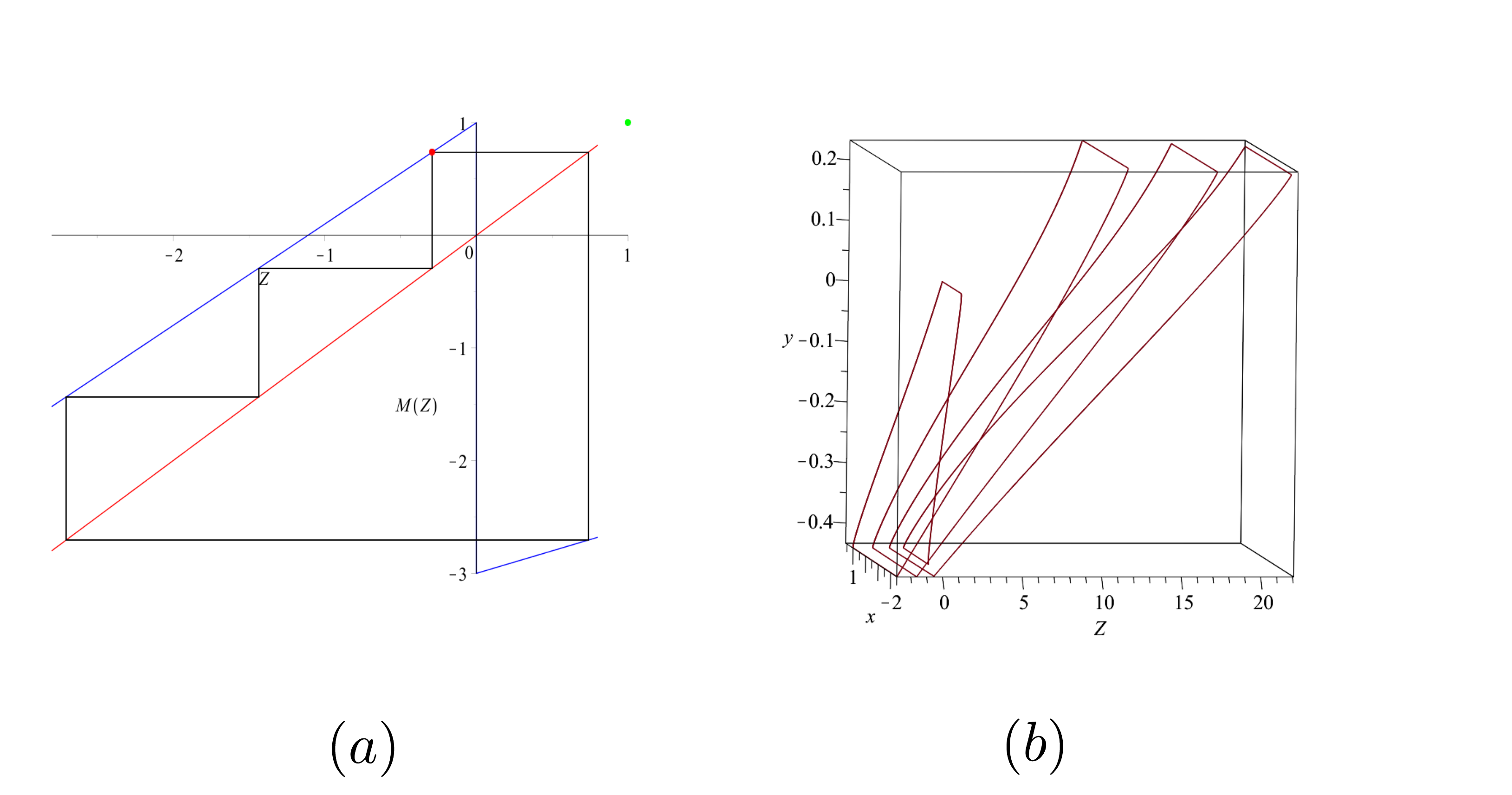} 
\caption{MMO of signature $3^1$: (a) piecewise affine map and (b) associated three-dimensional slow-fast system.}
\label{fig:MMO3-1}
\end{figure}

\begin{figure}[H]
\includegraphics[width=12cm,center]{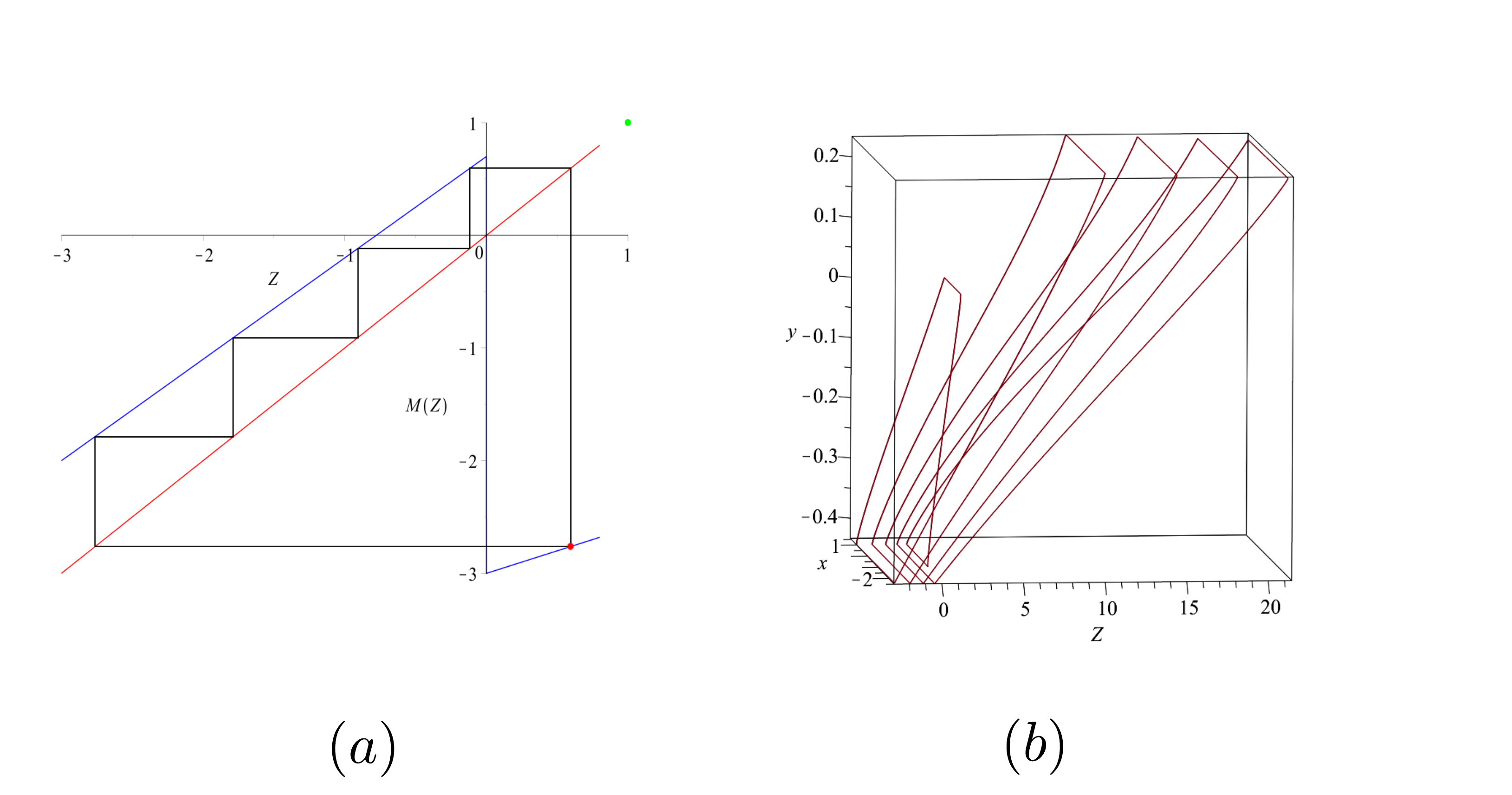} 
\caption{MMO of signature $4^1$: (a) piecewise affine map and (b) associated three-dimensional slow-fast system.}
\label{fig:MMO4-1}
\end{figure}

\begin{figure}[H]
\includegraphics[width=12cm,center]{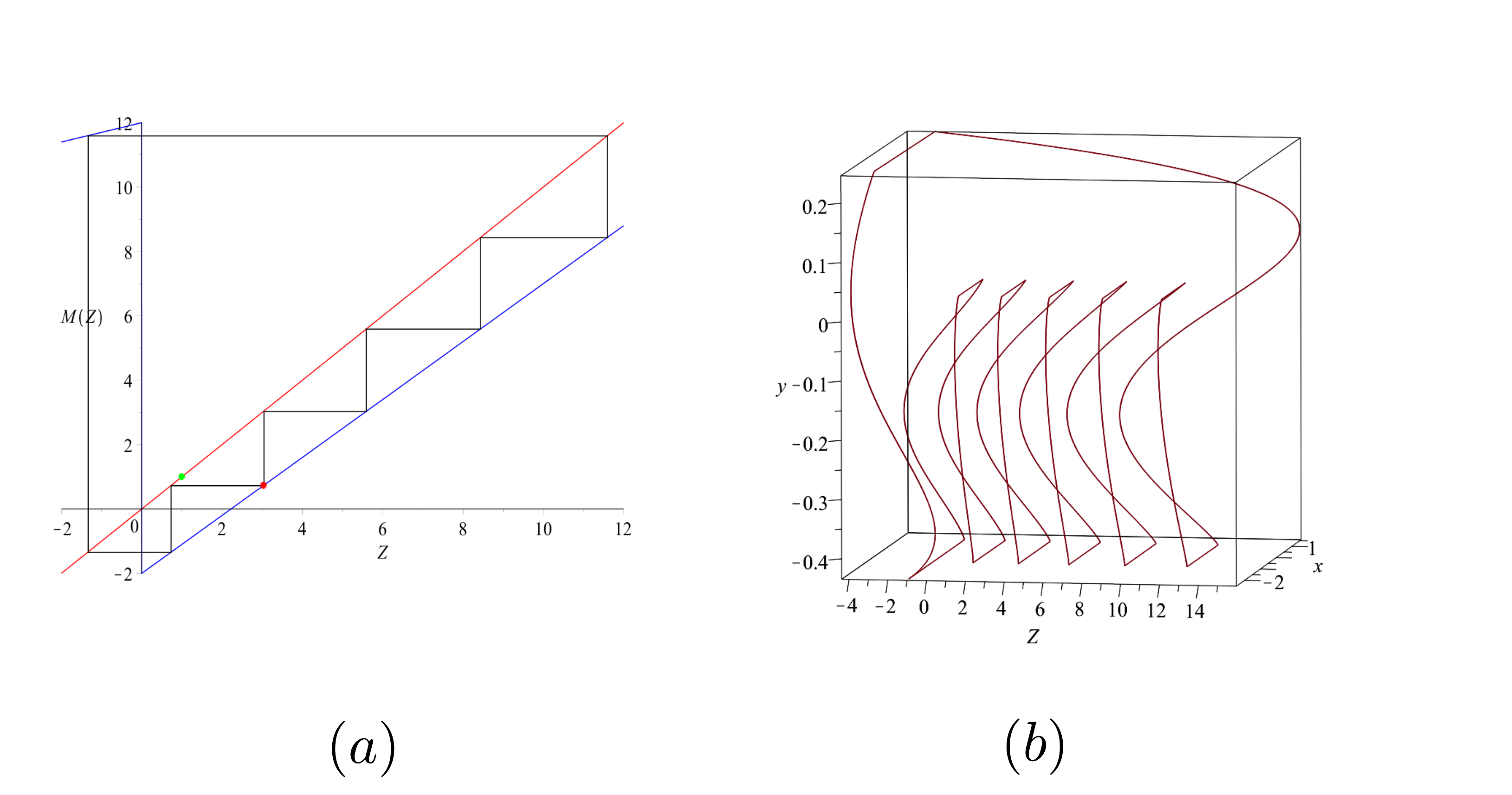} 
\caption{MMO of signature $1^5$: (a) piecewise affine map and (b) associated three-dimensional slow-fast system.}
\label{fig:MMO1-5}
\end{figure}

\begin{figure}[H]
\includegraphics[width=14cm,center]{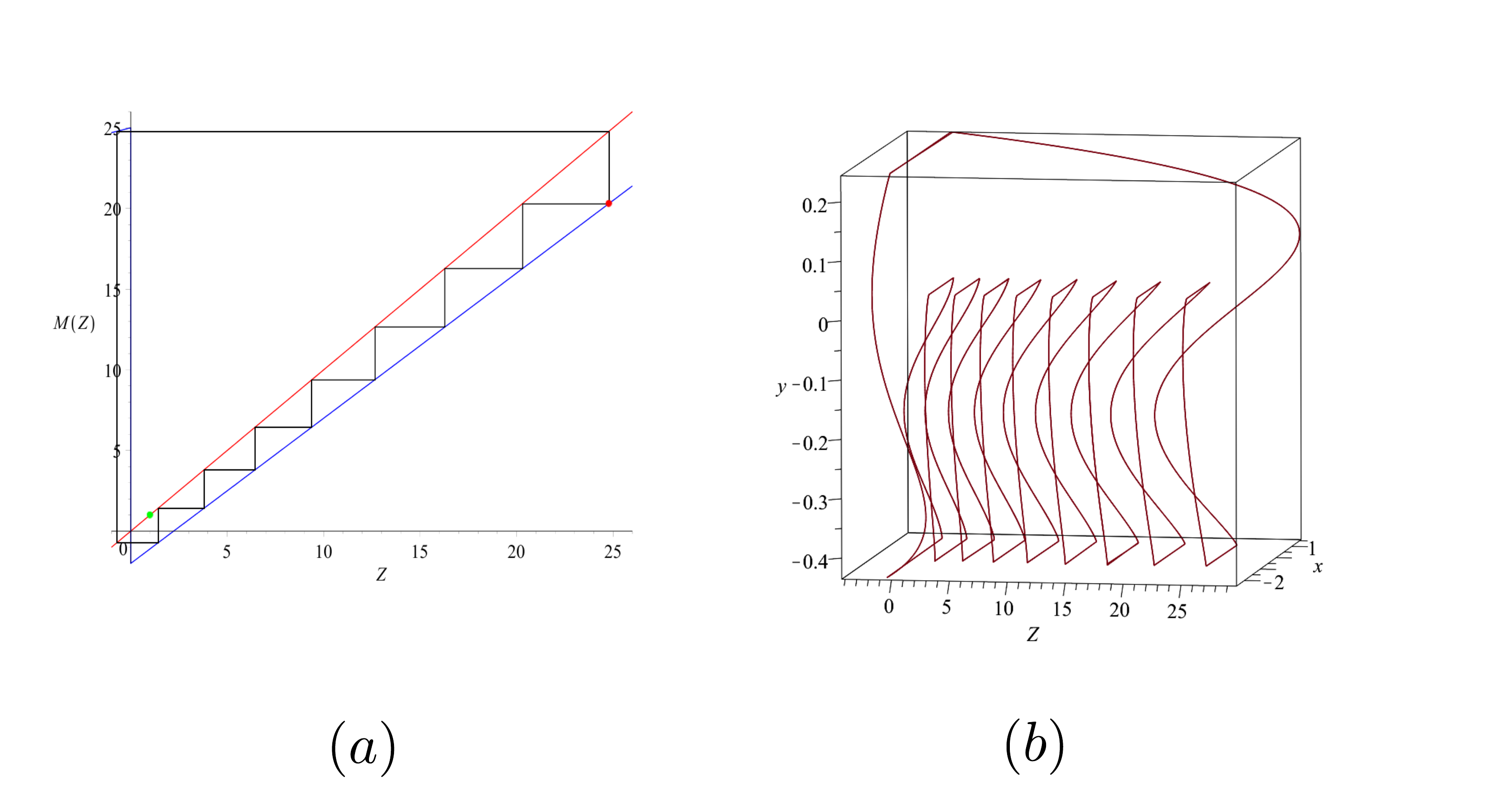} 
\caption{MMO of signature $1^8$: (a) piecewise affine map and (b) associated three-dimensional slow-fast system.}
\label{fig:MMO1_8}
\end{figure}

\subsection{The At Most \& At Least Lemma}\label{sec:atmost+atleast}

In this section, we give conditions on the coefficients in the definition of the PAM $M$ in \eqref{eq:PAM_prototype} that guarantee the occurrence of certain numbers of
LAOs ($Z<0$) or SAOs ($Z>0$) in a periodic MMO generated by the PAM $M$ in
\eqref{eq:PAM_prototype}. To that end, we apply results of \cite{RPB12}; in a first step, we
transform $M$ into the form considered there, via $a = a_{11}$, $b = a_{21}$, $\mu = a_{12}$, and $l= a_{22}-a_{12}$:
\begin{equation}
M(Z)=\begin{cases}
aZ + \mu & \quad\text{for }Z<0, \\
bZ + \mu + l & \quad\text{for }Z>0.
\end{cases}
\label{trans_PAM}
\end{equation}
To ensure the stability of the resulting MMOs, we assume that $a$ and $b$ take values in the interval $(0,1)$. Then, the parameter $l$ represents the height of the jump at $Z=0$, while the parameter $\mu$ will be varied. As explained in \cite{RPB12}, we restrict to $0< \mu < -l$, in which case \eqref{trans_PAM} has no fixed points and periodic orbits are possible.

 The following result then gives conditions on the control parameter $\mu$ for at most, or at least, $L$ consecutive LAOs, respectively $s$ SAOs, to appear in a periodic MMO for $M$. 
\begin{proposition}[At Most \& At Least Lemma \cite{RPB12}] Let $M$ be as defined in \eqref{trans_PAM}. Then, the following statements hold true.
	\begin{itemize}
\item[(1)]	When
		$\mu \leq  \frac{	-la^{L-1}	}{	a^{L-1}b + \sum^{L-1}_{k=0}a^{k}		}=:\mu_1$,
		then at least $L$ consecutive LAOs appear in a periodic MMO of $M$.
		When
		$\mu > \frac{	-la^{L}	}{	\sum^{L}_{k=0}a^{k}	}=:\mu_2$,
		then at most $L$ consecutive LAOs appear in a periodic MMO of $M$.
\item[(2)]	When
		$\mu < \frac{-l \sum^{s-1}_{k=0} b^k}{ \sum^{s}_{k=0} b^k }$,
		then at most $s$ consecutive SAOs appear in a periodic MMO of $M$.
		When 
                $\mu\ge \frac{-l\big[\sum^{s-1}_{k=0}b^{k}+b^{s-1}(a-1)\big]	}{	b^{s-1}a + \sum^{s-1}_{k=0}b^{k}}$,
		then at least $s$ consecutive SAOs appear in a periodic MMO of $M$.
		\label{atmost+atleast}
	\end{itemize}
\end{proposition}
Given Proposition~\ref{atmost+atleast}, it can be shown \cite{RPB12} that for $\mu\in(\mu_2,\mu_1]$, 
the only possible periodic MMO for $M$ is the one with signature $L^1$. Similarly, we can determine intervals for $\mu$ on which periodic MMOs with signature $1^s$ exist. We summarise a sample of MMO signatures, and the corresponding parameter regimes, in Table~\ref{table:param_int} below. Here, the relevant $\mu$-intervals are obtained from Proposition~\ref{atmost+atleast}; throughout, we find agreement between the theory (``Predicted $\mu$") and our numerics (``Actual $\mu$").

{\renewcommand{\arraystretch}{1.2}
	\begin{table}[h!]
	\setlength{\tabcolsep}{2pt}
		\centering 
		\begin{tabular}{|c|c|c|c|c|c|c|c|c|c|}
			\hline
			\hline                        
\tiny{\textbf{Signature}} &$a$ &$b$ &$l$ &\tiny{\textbf{Predicted $\mu$}}	&\tiny{\textbf{Actual $\mu$}}	&\textbf{$\alpha$}	&\textbf{$\beta$}	&\textbf{$\kappa$}	&\textbf{$\lambda$}		\\ [0.5ex]
			\hline
			\hline
			
	 	$1^2$	&0.3 	&0.9		&-5		&$\mu \in \big(2.9262,3.5055	\big]$	&$3$		&0.8743&0.0241&28.23&73.18			\\[2pt]
\hline	$1^3$	&0.3 	&0.9		&-9		&$\mu \in \big(6.5313,7.0921\big]$	&$7$		&0.8743&0.0241&30.1744&90.4070		\\[2pt]
\hline	$1^4$	&0.3    &0.9		&-11		&$\mu \in \big(8.8076,9.2376\big]$	&$9$		&0.8743&0.0241&31.1434&-99.0172		\\[2pt]   
\hline 	$1^8$	&0.3		&0.9		&-2.28	&$\mu \in (2.0932,2.1197]$			&$2.1$	&0.8743&0.0241&3.4279&-14.4651		\\[2pt]
\hline	$1^9$	&0.3		&0.9		&-2.68	&$\mu \in (2.4955,2.5205]$			&$2.5$	&0.8743&0.0241&3.6217&-16.1861		\\[2pt]
\hline	$1^{25}$	&0.5		&0.9	4	&-15.25	&$\mu \in \big(14.9889,15.0064\big]$	&$15$&0.5025&0.0152&15.8532&-177.4797		\\[2pt]
\hline  $2^1$	&0.9		&0.8		&-7.2	&$\mu \in \big[2.1520,2.4732\big)$	&$2.2$	&-0.0610&0.2430&24.4916&-96.1819		\\[2pt]
\hline	$3^1$	&0.9		&0.8		&-6.5	&$\mu \in \big[1.3778,1.5678\big)$	&$1.5$	&-0.0610&0.2430&24.5673&-81.8569		\\[2pt]
\hline	$6^1$	&0.9		&0.8		&-5.6	&$\mu \in \big[0.5704,0.6410\big)$	&$0.6$	&-0.0610&0.2430&24.6646&-63.4391		\\[2pt]  
\hline	$8^1$	&0.9		&0.9		&-6.5	&$\mu \in [0.4675,0.5075)$			&$0.5$	&0.0147&0.1104&65.5190&-462.9354		\\[2pt] 
\hline	$9^1$	&0.9		&0.9		&-9.6	&$\mu \in [0.5710,0.6344)$			&$0.6$	&0.0147&0.1104&98.1512&-683.7200		\\[2pt]
    \hline			
    \end{tabular}
    \caption{Signatures of the form $L^1$ and $1^s$ generated by (\ref{VF_represent}) and the corresponding $\mu$-intervals, as determined from Proposition~\ref{atmost+atleast}.}
    \label{table:param_int}
	\end{table}
}
Given the above, it is natural to ask whether MMOs with signature $L^s$ for $L>1$
and $s>1$ can be found in the present context. Following again \cite{RPB12}, it can be shown 
that stable periodic MMOs with such signatures cannot occur; we outline the argument here for completeness.
 In \cite{RPB12}, an orbit $\mathcal{O}$ is called \textit{admissible} if the $\mu$-interval for which $\mathcal{O}$ exists is non-empty. Then, their Lemma 2 states that ``\textit{for any admissible orbit $\mathcal{O}$, its pattern cannot contain consecutive $L$s and consecutive $R$s simultaneously}", where $L$ and $R$ denote numbers of LAOs and
SAOs in $\mathcal{O}$, respectively. The proof of Lemma~2 is by contradiction: if one assumes that an orbit with signature $L^R$ is actually possible, one concludes that, necessarily, $a,b>1$ in
Equation~\eqref{trans_PAM}; however, that contradicts the underlying assumption of $a,b \in (0,1)$ which is imposed in \cite{RPB12}.

In fact, since $a=a_{11}$ and $b=a_{21}$, $a,b>1$ would also imply instability of the corresponding MMO in \eqref{trans_PAM}; recall the proof of Theorem~\ref{association_thm} in Section~\ref{proof of assosciation theorem}. Hence, it is natural to assume that $0<a,b<1$ in our own analysis, as well, in which case the existence of more ``exotic" stable periodic MMOs with general signature $L^s$ can be ruled out. 

\subsection{Crossover signatures}\label{sec:crossover}

Given our numerical results in the previous two subsections, it is natural to ask what happens between two ``consecutive" signatures, i.e., how the shape of an MMO changes as orbits cross over from a cycle of signature $L^{s}$ to one of signature $L^{s+1}$ or, equivalently, from one of signature  $L^{s}$ to one of signature $(L+1)^{s}$.
Motivated again by results of \cite{RPB12} -- see, in particular, Lemma 4, Figure 3, and Note 2 therein -- we observe the existence of so-called ``crossover signatures" inside ``intermediate neighbourhoods" for some of the corresponding parameters in the definition of the transformed PAM $M$ in \eqref{trans_PAM}. (In \cite{RPB12}, the existence of similar regions, named ``molecular regions" there, is concluded.) These observations lead to the conclusion that MMO signatures are not arranged in a monotonous way, as far as the number of LAOs or SAOs therein is concerned. For illustration, we showcase a simple case here, namely, an MMO of signature $1^4$, which can be obtained from the following PAM,
\begin{equation}
M_{1^{4}}(Z)=\begin{cases}
0.9 Z + 6 & \quad\text{for }Z<0,	\\
0.85Z - 1 & \quad\text{for }Z>0,
\end{cases}
\end{equation}
with corresponding parameter values $\alpha \approx -0.0220$, $\beta \approx 0.1747$, $\kappa \approx 7.7321$, and $\lambda \approx -233.3068$ in the associated slow-fast vector field that is determined by \eqref{VF_represent}.
It is straightforward to obtain an MMO with the ``consecutive signature", namely $1^5$, in the following PAM:
\begin{equation}
M_{1^{5}}(Z)  =	\begin{cases}
0.9 Z + 7.2 & \quad\text{for }Z<0, \\
0.85Z - 1	& \quad\text{for }Z>0,
\end{cases}
\end{equation}
with parameter values $\alpha \approx -0.0220$, $\beta \approx 0.1747$, $\kappa \approx 7.9239$, and $\lambda \approx -275.3021$.

\begin{figure}[H]
\includegraphics[width=12cm,center]{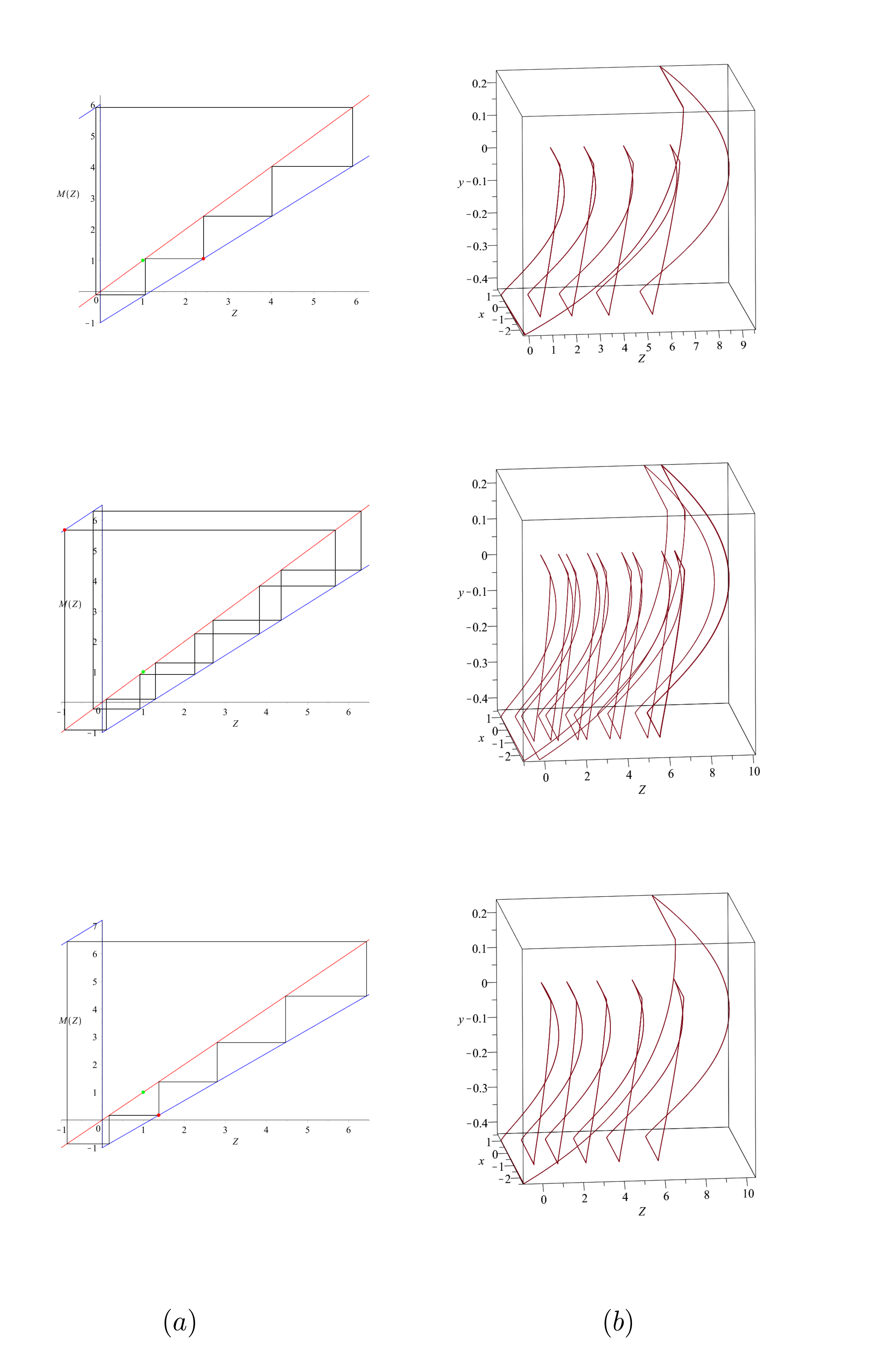} 
\caption{Transition from MMOs with signature $1^4$ to $1^5$ via the "crossover" signature $1^41^5$ in (a) the PAM defined in \eqref{eq:PAM_prototype} and (b) its associated three-dimensional slow-fast system. The signature $1^41^5$ was found for $\kappa = 7.8120$ and $\lambda = -250.8049$, with $\mu=6.5$.}
\label{fig:crossover1415}
\end{figure}

Noting that the numerical values of the parameters $\alpha$ and $\beta$ that determine $M_{1^{4}}$
and $M_{1^{5}}$ are almost identical while $\kappa$ and $\lambda$ vary, we take
$\kappa \in (7.7321, 7.9239)$ and $\lambda \in (-275.3021, -233.3068)$, which generates
the ``crossover signature" $1^{4}1^{5}$ for $\mu=6.5$, as shown in Figure~\ref{fig:crossover1415}(a).

Here, it is important to emphasise that these mixed signatures do not contradict the At Most \& At Least Lemma, Proposition~\ref{atmost+atleast}. Rather, for a fixed choice of the pivotal quantities 
$\alpha$, $\beta$, $\kappa$, and $\lambda$, we obtain a hierarchy of disjoint $\mu$-intervals that
correspond to mixed-mode signatures of the form $L^s$ from Proposition~\ref{atmost+atleast}. 
``Crossover" signatures are found for $\mu$ chosen in the complements of those intervals; from a practical point of view, our choice of the pivotal quantities is guided by where the adjacent, ``simple'' signatures are found, whereupon $\mu$ can be fixed from the At Most \& At Least Lemma. See
Table~\ref{table:param_int} for a specification of the corresponding $\mu$-intervals.

\begin{figure}[H]
\includegraphics[width=12cm,center]{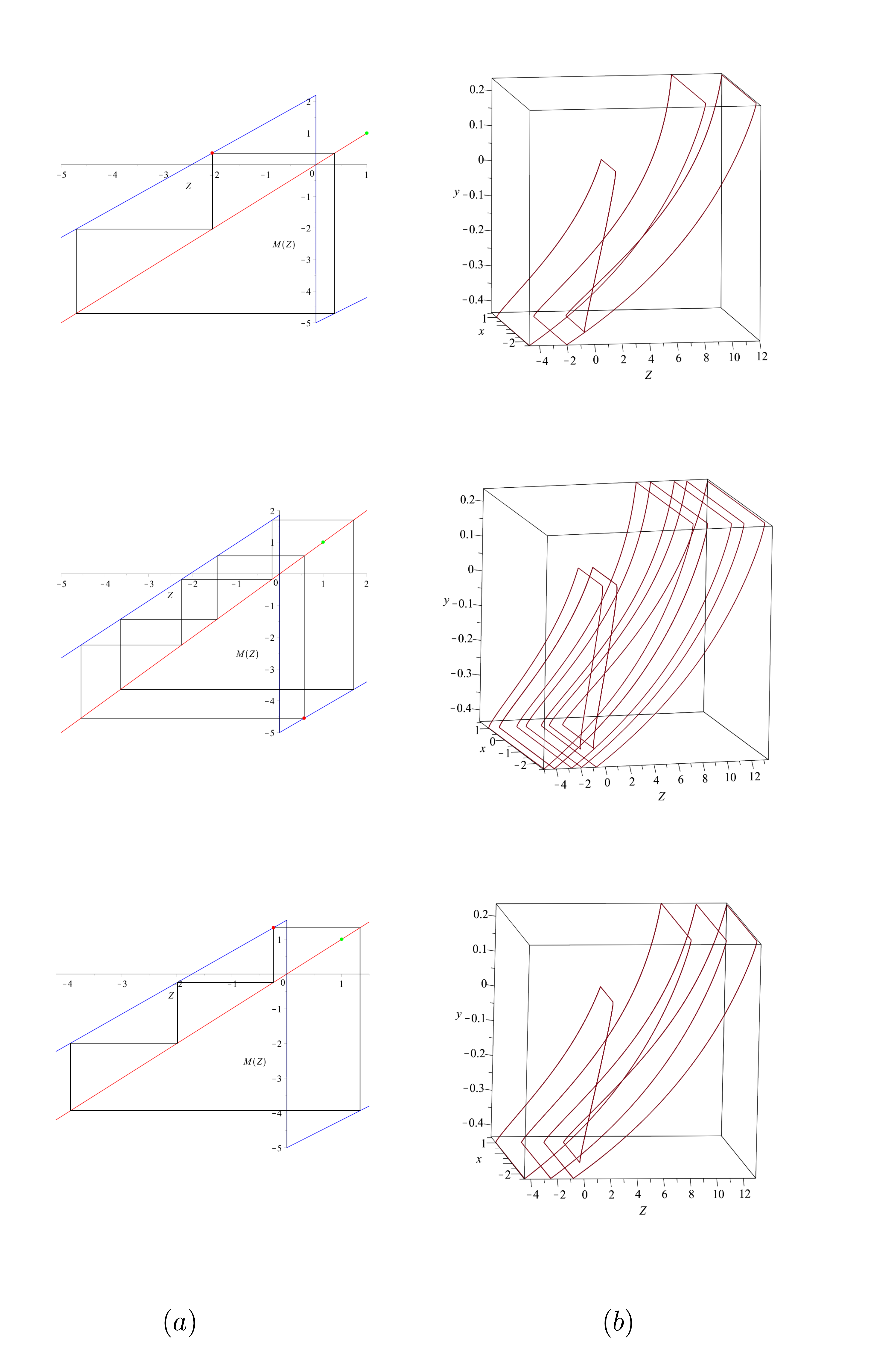} 
\caption{Transition from MMOs with signature $2^1$ to $3^1$ via the "crossover" signature $2^13^1$ in (a) the PAM defined in \eqref{eq:PAM_prototype} and (b) its associated three-dimensional slow-fast system. The signature $2^13^1$ was found for $\alpha = -0.0610$, $\beta = 0.2430$, $\kappa = 24.5348$, and $\lambda = -87.9962$, with $\mu=1.8$.}
\label{fig:crossover2131}
\end{figure}

Following the same procedure as above, we were able to detect intermediate neighbourhoods for the signature $2^1$ crossing over to $3^1$; see Figure~\ref{fig:crossover2131} for an illustration.

Again, we first consider a PAM which realises the signature $2^1$:
\begin{equation}
M_{2^1}(Z)=\begin{cases}
0.9 Z + 2.2 & \quad\text{for }Z<0, \\
0.8 Z - 5 & \quad\text{for }Z>0,
\end{cases}
\end{equation}
where $\alpha = -0.0610$, $\beta = 0.2430$, $\kappa = 24.4916$, and $\lambda = -96.1819$, as well as
a map which generates the ``consecutive" signature $3^1$:
\begin{equation}
M_{3^1}(Z)=\begin{cases}
0.9 Z + 1.5 & \quad\text{for }Z<0, \\
0.8 Z - 5 & \quad\text{for }Z>0,
\end{cases}
\end{equation}
with $\alpha = -0.0610$, $\beta = 0.2430$, $\kappa = 24.5673$, and $\lambda = -81.8569$.

 Naturally, here too the numerical values of the parameters $\alpha$ and $\beta$ that determine $M_{2^1}$ and $M_{3^1}$ are almost indistinguishable. If we then pick $\kappa\in (24.4916, 24.5673)$ and $\lambda\in (-96.1819, -81.8569)$, we observe an MMO with crossover signature $2^{1}3^{1}$ for $\mu=1.8$, as shown in Figure~\ref{fig:crossover2131}.

Note that for both examples, we took $a_{11}=a<1$ and $a_{21}=b<1$, in accordance with \cite{RPB12}; a more general choice of coefficients, with $a^{L}b^s<1$, still yields stable mixed-mode orbits for the PAM $M$ in \eqref{trans_PAM}. However, as stated above, that scenario is excluded in \cite{RPB12} and hence cannot be considered within the framework of the At Most \& At Least Lemma, Proposition~\ref{atmost+atleast}.

\begin{remark}\label{rem:crossover}
Numerical evidence suggests that we only encounter a combination of either one of the adjacent signatures between two ``simple", consecutive signatures: for example, picking $1^s$ and $1^{s+1}$, in the ''intermediate neighbourhoods" we can only expect MMOs of signature $(1^{s} )^{i} (1^{s+1})^{j}$, with finite multiplicity $0<i,j < \infty$. This last assertion is easier to observe for small $s$, since the parameter intervals corresponding to such ``crossover" signatures tend to shrink with increasing $s$. An analogous assertion applies to MMOs with signature $L^1$. 
\end{remark}

\section{Conclusions}\label{Conclusions}

In this paper, we have introduced a novel ``jump mechanism" for the generation of mixed-mode oscillations (MMOs) in a family of three-dimensional slow-fast systems of singular perturbation type. 
In marked contrast to the canard-based mechanism that is typically invoked to explain mixed-mode dynamics in such systems, we do not assume the presence of a folded singularity on any of the fold lines in our system at which normal hyperbolicity is lost; 
in fact, we require all such fold lines to consist of jump points only. 
Correspondingly, the SAO component in the resulting mixed-mode trajectories is then also of relaxation type, with an amplitude that is $O(1)$ in the singular perturbation parameter $\epsilon$.
At this point, we remark that it is possible to obtain quantitative information on the MMOs constructed here. In particular, the amplitudes of both the LAO and SAO components can be determined by observing the height of the fold lines $L^1$ and $L^3$ with respect to the $y$-coordinate. The corresponding periods can be approximated from the transition times on normally attracting portions of the critical manifold in a manner similar to that of Section ~\ref{sec:derivePAM}; given a mixed-mode trajectory of signature $L^s$, the overall period would be found by multiplying the periods of one LAO and one SAO with $L$ and $s$, respectively, before adding them. The details are left to the interested reader.

As our principal result, we have established a two-way correspondence between our family of slow-fast systems and a class of one-dimensional piecewise affine maps (PAMs) which are naturally associated to each other. 
In particular, we have shown that for every such PAM that exhibits an MMO with a certain given signature, there exists a slow-fast system that can be associated to it and {\it vice versa}, given that certain conditions are met. 
Thus, we have reduced the study of MMOs in a relatively broad family of singularly perturbed three-dimensional systems to the well-developed theory of one-dimensional maps.
We were able to verify our own results numerically, showing that they are consistent with those obtained in \cite{RPB12} -- and, in particular, with the At Most \& At Least Lemma -- in the process. 

Naturally, a number of questions arise from the present analysis. 
The first of these concerns an in-depth investigation of a neighbourhood of the singular orbit $\Gamma_0$ defined in Assumption~\ref{A5}, as well as of the corresponding discontinuity in the associated PAM, where canard phenomena could occur. We conjecture that this discontinuity gives rise to a canard explosion which determines the interchange between LAOs and SAOs in the resulting mixed-mode time series.

Next, it seems natural to comment on the interplay between $\epsilon$ and $\delta$ in 
our prototypical family of slow-fast systems in Equation~\eqref{VFfast}, which we restate below for
reference:
\begin{align*}
x' &=y-F(x,z,\epsilon,\delta)=:f(x,y,z,\epsilon,\delta), \\
y' &=\epsilon g_1(x,y,z,\epsilon,\delta), \\
z' &= \epsilon g_2(x,y,z,\epsilon,\delta).
\end{align*}
The form of the above two-parameter singular perturbation problem 
implies that our choice of $\delta$ will mostly affect the slow flow along normally hyperbolic (attracting) portions of the corresponding critical manifold $\mathcal{S}$,
away from the fold lines $L_i$. Conversely, it also appears that any restrictions on the magnitude of $\epsilon$ are only due to the jump behaviour at those lines; both observations are corroborated by
numerical experimentation in {\sc Maple}.

Now, the initial rescaling of $z$ with respect to $\delta$, as illustrated in Section~\ref{sec:derivePAM}, implies that the relevant $z$-window for our analysis is $\mathcal{O}(\delta)$ wide.
While $Z<0$ — in the rescaled $Z$-variable — results in LAOs in the corresponding mixed-mode time series, whereas $Z>0$ yields SAOs, that classification is only true in the singular limit, i.e. for $\epsilon = 0$. To specify the interplay between $\delta$ and $\epsilon$ away from that limit, one
would need to ``blow up” (desingularise) the flow of \eqref{VFfast} in the vicinity of the degenerate
point $P_c$; as $\epsilon$ would receive triple the weight of
$\delta$ in the corresponding blow-up transformation, it would then follow that we not only have to avoid $P_c$ itself, but also an $\mathcal{O}(\epsilon^{1/3})$ ``hole" around that
point, in order not to have to consider canard phenomena in our analysis. When $\epsilon^{1/3} > \delta$, on the other hand, the Poincar\'e map $\Pi$ associated to \eqref{VFfast} is likely to return the flow inside this \textit{canard hole}. 
Correspondingly, in Section~\ref{sec:numerics}, we fixed $\epsilon = 10^{-7}$ in our numerics, as we had experimentally concluded that the optimal choice of $\delta$ is somewhere in the region of $5 \times 10^{-3}$, in agreement with the above reasoning.

The next question that comes to mind, which concerns the patterns that the resulting MMO signatures follow, is motivated by results of Freire and Gallas in \cite{FG11}. There, it is shown that the number of SAOs in a given mixed-mode orbit is not arbitrary, but that it is organised in a pattern dictated by a so-called Stern-Brocot tree. It would seem natural to investigate whether similar number-theoretical
arguments can be applied in the context of the family of slow-fast systems studied in the present work. Preliminary analysis seems to suggest that the ``crossover” signatures observed between
simple patterns of the form $1^s$ or $L^1$ are relatively regular; recall Remark~\ref{rem:crossover}.

Our final remark concerns the prototypical family of slow-fast systems in \eqref{VFfast}, and specifically those which have the property in \eqref{eq:g2}. We are confident that, near the fold lines $L_{1}$, $L_{3}$, and $L_{4}$, a suitable Fenichel-like normal form can be derived. Should that expectation be verified, it would follow that one could apply our approach to connect the well-developed theory of one-dimensional PAMs with the vast family of slow-fast systems that can be brought into said normal form.

\section*{Acknowledgements}
The authors thank the School of Mathematics at the University of Edinburgh for its hospitality during several research visits. In particular, we are grateful to Panagiotis Kaklamanos for his fruitful and meticulous comments, as well as to the entire Edinburgh Dynamical Systems Study Group for general feedback on a draft version of the paper.


\bibliographystyle{abbrv}
\bibliography{references}

\begin{thebibliography}{10}

\bibitem{BKYY2000}
S.~Banerjee, M.~S. Karthik, G.~Yuan, and J.~A. Yorke.
\newblock Bifurcations in one-dimensional piecewise smooth maps---theory and
  applications in switching circuits.
\newblock {\em IEEE Trans. Circuits Systems I Fund. Theory Appl.},
  47(3):389--394, 2000.

\bibitem{BENOIT}
E.~Beno\^{i}t.
\newblock Syst\`emes lents-rapides dans {${\mathbb R}^{3}$} et leurs canards.
\newblock In {\em Third {S}chnepfenried geometry conference, {V}ol. 2
  ({S}chnepfenried, 1982)}, volume 109 of {\em Ast\'{e}risque}, pages 159--191.
  Soc. Math. France, Paris, 1983.

\bibitem{BCDD81}
E.~Benoit, J.~Louis~Callot, F.~Diener, and M.~Diener.
\newblock Chasse au canard.
\newblock {\em Collectanea Mathematica}, 32, 01 1981.

\bibitem{BKW06}
M.~Brøns, M.~Krupa, and M.~Wechselberger.
\newblock Mixed mode oscillations due to the generalized canard phenomenon.
\newblock {\em Fields Institute Communications}, 49:39--63, 10 2006.

\bibitem{DGKKOW12}
M.~Desroches, J.~Guckenheimer, B.~Krauskopf, C.~Kuehn, H.~M. Osinga, and
  M.~Wechselberger.
\newblock Mixed-mode oscillations with multiple time scales.
\newblock {\em SIAM Rev.}, 54(2):211--288, May 2012.

\bibitem{DIBERNARDO}
M.~{di Bernardo}, C.~Budd, A.~Champneys, and P.~Kowalczyk.
\newblock {\em Piecewise-smooth dynamical systems: theory and applications}.
\newblock Applied Mathematical Sciences. Springer, 2008.

\bibitem{D84}
M.~Diener.
\newblock The canard unchained or how fast/slow dynamical systems bifurcate.
\newblock {\em Math. Intelligencer}, 6(3):38--49, 1984.

\bibitem{D94}
M.~Diener.
\newblock Regularizing microscopes and rivers.
\newblock {\em SIAM J. Math. Anal.}, 25(1):148--173, 1994.

\bibitem{DR96}
F.~Dumortier and R.~Roussarie.
\newblock Canard cycles and center manifolds.
\newblock {\em Mem. Amer. Math. Soc.}, 121(577):x+100, 1996.
\newblock With an appendix by Cheng Zhi Li.

\bibitem{DRBA08}
P.~S. Dutta, B.~Routroy, S.~Banerjee, and S.~S. Alam.
\newblock On the existence of low-period orbits in {$n$}-dimensional piecewise
  linear discontinuous maps.
\newblock {\em Nonlinear Dynam.}, 53(4):369--380, 2008.

\bibitem{F72}
N.~Fenichel.
\newblock Persistence and smoothness of invariant manifolds for flows.
\newblock {\em Indiana Univ. Math. J.}, 21:193--226, 1972.

\bibitem{F79}
N.~Fenichel.
\newblock Geometric singular perturbation theory for ordinary differential
  equations.
\newblock {\em J. Differential Equations}, 31(1):53--98, 1979.

\bibitem{FG11}
J.~G. Freire and J.~A. Gallas.
\newblock Stern–{B}rocot trees in cascades of mixed-mode oscillations and
  canards in the extended {B}onhoeffer–van der {P}ol and the
  {F}itz{H}ugh–{N}agumo models of excitable systems.
\newblock {\em Physics Letters A}, 375(7):1097 -- 1103, 2011.

\bibitem{H00}
T.~Hayashi.
\newblock Mixed-mode oscillations and chaos in a glow discharge.
\newblock {\em Phys. Rev. Lett.}, 84:3334--3337, Apr 2000.

\bibitem{IMMTDZ11}
C.~Iglesias, C.~Meunier, M.~Manuel, Y.~Timofeeva, N.~Delestr{\'e}e, and
  D.~Zytnicki.
\newblock Mixed mode oscillations in mouse spinal motoneurons arise from a low
  excitability state.
\newblock {\em Journal of Neuroscience}, 31(15):5829--5840, 2011.

\bibitem{IK2020}
N.~Inaba and T.~Kousaka.
\newblock Nested mixed-mode oscillations.
\newblock {\em Physica D: Nonlinear Phenomena}, 401:132152, 2020.

\bibitem{JMBNR13}
N.~D. Jimenez, S.~Mihalas, R.~Brown, E.~Niebur, and J.~Rubin.
\newblock {Locally contractive dynamics in generalized integrate-and-fire
  neurons}.
\newblock {\em SIAM J Appl Dyn Syst}, 12(3):1474--1514, Sep 2013.

\bibitem{JONES}
C.~K. R.~T. Jones.
\newblock {\em Geometric singular perturbation theory}, pages 44--118.
\newblock Springer Berlin Heidelberg, Berlin, Heidelberg, 1995.

\bibitem{chemists}
K.~Kovacs, M.~Leda, V.~K. Vanag, and I.~R. Epstein.
\newblock Small-amplitude and mixed-mode p{H} oscillations in the
  bromate−sulfite−ferrocyanide−aluminum(iii) system.
\newblock {\em The Journal of Physical Chemistry A}, 113(1):146--156, 2009.
\newblock PMID: 19086810.

\bibitem{KPK08}
M.~Krupa, N.~Popović, and N.~Kopell.
\newblock Mixed-mode oscillations in three time-scale systems: {A} prototypical
  example.
\newblock {\em SIAM J. Applied Dynamical Systems}, 7:361--420, 01 2008.

\bibitem{KS01}
M.~Krupa and P.~Szmolyan.
\newblock Extending geometric singular perturbation theory to nonhyperbolic
  points---fold and canard points in two dimensions.
\newblock {\em SIAM J. Math. Anal.}, 33(2):286--314, 2001.

\bibitem{MN09}
S.~Mihalaş and E.~Niebur.
\newblock A generalized linear integrate-and-fire neural model produces diverse
  spiking behaviors.
\newblock {\em Neural Computation}, 21(3):704--718, 2009.
\newblock PMID: 18928368.

\bibitem{MSLG97}
A.~Milik, P.~Szmolyan, H.~Löffelmann, and E.~Gröller.
\newblock Geometry of mixed-mode oscillations in the 3-{D} autocatalator.
\newblock {\em International Journal of Bifurcation and Chaos}, 8, 01 1997.

\bibitem{PRV16}
R.~Prohens, A.~Teruel, and C.~Vich.
\newblock Slow--fast $n$-dimensional piecewise linear differential systems.
\newblock {\em Journal of Differential Equations}, 260(2):1865--1892, 2016.

\bibitem{RPB12}
B.~Rajpathak, H.~K. Pillai, and S.~Bandyopadhyay.
\newblock Analysis of stable periodic orbits in the one dimensional linear
  piecewise-smooth discontinuous map.
\newblock {\em Chaos: An Interdisciplinary Journal of Nonlinear Science},
  22(3):033126, 2012.

\bibitem{SSI15}
K.~Shimizu, M.~Sekikawa, and N.~Inaba.
\newblock Experimental study of complex mixed-mode oscillations generated in a
  {B}onhoeffer-van der {P}ol oscillator under weak periodic perturbation.
\newblock {\em Chaos (Woodbury, N.Y.)}, 25:023105, 02 2015.

\bibitem{SW04}
P.~Szmolyan and M.~Wechselberger.
\newblock Relaxation oscillations in {${\mathbb R}^3$}.
\newblock {\em J. Differential Equations}, 200(1):69--104, 2004.

\bibitem{GKE17}
B.~V-Ghaffari, M.~Kouhnavard, and S.~M. Elbasiouny.
\newblock Mixed-mode oscillations in pyramidal neurons under antiepileptic drug
  conditions.
\newblock {\em PLOS ONE}, 12(6):1--20, 06 2017.

\bibitem{VL97}
A.~Venkatesan and M.~Lakshmanan.
\newblock Bifurcation and chaos in the double-well {D}uffing--van der {P}ol
  oscillator: {N}umerical and analytical studies.
\newblock {\em Phys. Rev. E}, 56:6321--6330, Dec 1997.

\bibitem{W05}
M.~Wechselberger.
\newblock Existence and bifurcation of canards in {${\mathbb R}^3$} in the case
  of a folded node.
\newblock {\em SIAM Journal on Applied Dynamical Systems}, 4(1):101--139, 2005.

\bibitem{YO96}
T.~Yamaguchi and H.~Ohtagaki.
\newblock The order of appearance of oscillation modes of a piecewise linear
  map.
\newblock {\em J. Phys. Soc. Japan}, 65(11):3500--3512, 1996.

\end{thebibliography}
\end{document}